\documentclass[12pt]{article}

\usepackage{enumitem}

\usepackage{amsmath,amsthm,amssymb,amscd}
\usepackage[pdftex]{graphicx}
\usepackage{graphicx}
\usepackage{epsfig}
\usepackage{verbatim}
\usepackage{footmisc}
\usepackage{bbm,bm}
\usepackage{lscape}
\usepackage[longnamesfirst]{natbib}
\usepackage{color}
\usepackage{float}
\usepackage{longtable}
\usepackage{lipsum}

\newtheorem{Lemma}{Lemma}[section] 
\newtheorem{Theorem}{Theorem}[section] 

\makeatletter

\@addtoreset{equation}{section}
\makeatother

\renewcommand{\baselinestretch}{1.2}

\newcommand{\eps}{\varepsilon}

\def\bib{\par\noindent\hangindent=0.5 true cm\hangafter=1}

\newcommand{\sumi}{\sum_{i=1}^n}
\newcommand{\sumj}{\sum_{j=1}^n}

\newcommand{\R}{\mathbb{R}}

\newcommand{\bbeta}{\ensuremath{\boldsymbol{\beta}}}
\newcommand{\beps}{\ensuremath{\boldsymbol{\varepsilon}}}
\newcommand{\bx}{\bold{x}}
\newcommand{\bX}{\bold{X}}
\newcommand{\bY}{\bold{Y}}
\newcommand{\bz}{\bold{z}}
\newcommand{\bb}{\bold{b}}

\usepackage{color}

\begin{document}

\begin{titlepage}
\title{\bf Bootstrap of residual processes in regression: to smooth or not to smooth ?}
\author{{\large Natalie N\textsc{{eumeyer}}\footnote{Corresponding author; Department of Mathematics, University of Hamburg, Bundesstrasse 55, 20146 Hamburg, Germany, E-mail: neumeyer@math.uni-hamburg.de. Financial support by the DFG (Research Unit FOR 1735 {\it Structural Inference in
Statistics: Adaptation and Efficiency}) is gratefully acknowledged.}}
\and \addtocounter{footnote}{2}
 {\large Ingrid V\textsc{{an} K{eilegom}}\footnote{ORSTAT, KU Leuven, Naamsestraat 69, 3000 Leuven, Belgium, E-mail: ingrid.vankeilegom@kuleuven.be. Research financed by the European Research Council (2016- 2021, Horizon 2020 / ERC grant agreement No.\ 694409), and by IAP research network grant nr.\ P7/06 of the Belgian government (Belgian Science Policy).}}
\\[.5cm]}

\maketitle

\renewcommand{\baselinestretch}{1.1}

\begin{abstract}
In this paper we consider a location model of the form $Y = m(X) + \varepsilon$, where $m(\cdot)$ is the unknown regression function, the error $\varepsilon$ is independent of the $p$-dimensional covariate $X$ and $E(\varepsilon)=0$.  Given i.i.d.\ data $(X_1,Y_1),\ldots,(X_n,Y_n)$ and given an estimator $\hat m(\cdot)$ of the function $m(\cdot)$ (which can be parametric or nonparametric of nature), we estimate the distribution of the error term $\varepsilon$ by the empirical distribution of the residuals $Y_i-\hat m(X_i)$, $i=1,\ldots,n$.  To approximate the distribution of this estimator, Koul and Lahiri (1994) and Neumeyer (2008, 2009) proposed bootstrap procedures, based on smoothing the residuals either before or after drawing bootstrap samples.  So far it has been an open question whether a classical non-smooth residual bootstrap is asymptotically valid in this context.  In this paper we solve this open problem, and show that the non-smooth residual bootstrap is consistent.  We illustrate this theoretical result by means of simulations, that show the accuracy of this bootstrap procedure for various models, testing procedures and sample sizes.

\end{abstract}

\vspace*{.5cm}

\noindent{\bf Key words:} Bootstrap, empirical distribution function, kernel smoothing, linear regression, location model, nonparametric regression.

\end{titlepage}

\small
\normalsize
\addtocounter{page}{1}

\section{Introduction}  \label{introduction}

Consider the model
\begin{eqnarray} \label{model}
Y = m(X) + \varepsilon, 
\end{eqnarray}
where the response $Y$ is univariate, the covariate $X$ is of dimension $p \ge 1$, and the error term $\varepsilon$ is independent of $X$.  The regression function $m(\cdot)$ can be parametric (e.g.\ linear) or nonparametric of nature, and the distribution $F$ of $\eps$ is completely unknown, except that $E(\varepsilon)=0$.  The estimation of the distribution $F$ has been the object of many papers in the literature, starting with the seminal papers of Durbin (1973), Loynes (1980) and Koul (1987) in the case where $m(\cdot)$ is parametric, whereas the nonparametric case has been studied by Van Keilegom and Akritas (1999), Akritas and Van Keilegom (2001) and M\"uller et al.\ (2004), 
among others.   

The estimator of the error distribution has been shown to be very useful for testing hypotheses regarding several features of model (\ref{model}), like e.g.\ testing for the form of the regression function $m(\cdot)$ (Van Keilegom et al.\ (2008)), comparing regression curves (Pardo-Fern\'andez et al.\ (2007)), testing independence between $\eps$ and $X$ (Einmahl and Van Keilegom (2008) and Racine and Van Keilegom (2017)), testing for symmetry of the error distribution (Koul (2002), Neumeyer and Dette (2007)), 
among others. The idea in each of these papers is to compare an estimator of the error distribution obtained under the null hypothesis with an estimator that is not based on the null.  Since the asymptotic distribution of the estimator of $F$ has a complicated covariance structure, bootstrap procedures have been proposed to approximate the distribution of the estimator and the critical values of the tests.   

Koul and Lahiri (1994) proposed a residual bootstrap for linear regression models, where the bootstrap residuals are drawn from a smoothed empirical distribution of the residuals.   Neumeyer (2009) proposed a similar bootstrap procedure for nonparametric regression models.  The reason why a smooth bootstrap was proposed is that the methods of proof in both papers require a smooth distribution of the bootstrap error. 
Smooth residual bootstrap procedures have been applied by De Angelis et al.\ (1993), Mora (2005), Pardo-Fern\'andez et al.\ (2007), and Huskova and Meintanis (2009), among many others. 
 An alternative bootstrap procedure for nonparametric regression was proposed in Neumeyer (2008), where bootstrap residuals were drawn from the non-smoothed empirical distribution of the residuals, after which smoothing is applied on the empirical distribution of the bootstrap residuals. 
Further it has been shown that wild bootstrap in the context of residual-based procedures can only be applied for specific testing problems as testing for symmetry (Neumeyer et al.\ (2005), Neumeyer and Dette (2007)), whereas it is not valid in general (see Neumeyer (2006)).
 It has been an open question so far whether a classical non-smooth residual bootstrap is asymptotically valid in this context.   In this paper we solve this open problem, and show that the non-smooth residual bootstrap is consistent when applied to residual processes.   We will do this for the case of nonparametric regression with random design and for linear models with fixed design. Other models (nonparametric regression with fixed design, nonlinear or semiparametric regression,..) can be treated similarly. 
 The question whether smooth bootstrap procedures should be preferred over non-smooth bootstrap procedures has been discussed in different contexts, see Silverman and Young (1987) and Hall et al.\ (1989). 

The finite sample performance of the smooth and non-smooth residual bootstrap for residual processes has been studied by Neumeyer (2009).  The paper shows that for small sample sizes using the classical residual bootstrap version of the residual empirical process in the nonparametric regression context yields too small quantiles. However, as we will show in this paper, this problem is diminished for larger sample sizes and it is not very relevant when applied to testing problems. 

This paper is organized as follows.   In the next section, we show the consistency of the non-smooth residual bootstrap in the context of nonparametric regression.  In Section 3, the non-smooth residual bootstrap is shown to be valid in a linear regression model.  The finite sample performance of the proposed bootstrap is studied in a simulation study in Section 4.  All the proofs and the regularity conditions are collected in the Appendix.

\section{Nonparametric regression} \label{nonpreg}

We start with the case of nonparametric regression with random design.  The covariate is supposed to be one-dimensional.   To estimate the regression function we use a kernel estimator based on Nadaraya-Watson weights :
$$ \hat m(x) = \sumi \frac{k_h(x-X_i)}{\sumj k_h(x-X_j)} Y_i, $$
where $k$ is a kernel density function, $k_h(\cdot) = k(\cdot/h)/h$ and $h=h_n$ is a positive bandwidth sequence converging to zero when $n$ tends to infinity.   

Let $\hat \eps_i = Y_i - \hat m(X_i)$, $i=1,\ldots,n$.  In Theorem 1 in Akritas and Van Keilegom (2001) it is shown that the residual process $n^{-1/2} \sumi \big(I\{\hat\eps_i \le y\} - F(y)\big)$, $y \in \R$, converges weakly to a zero-mean Gaussian process $W(y)$ with covariance function given by 
\begin{eqnarray} \label{covW}
Cov\big(W(y_1),W(y_2)\big) &=& E \Big[\big(I\{\eps \le y_1\} + f(y_1) \eps\big) \, \big(I\{\eps \le y_2\} + f(y_2) \eps\big)\Big],
\end{eqnarray}
where $\eps$ has distribution function $F$ and density $f$. 

Neumeyer (2009) studied a smooth bootstrap procedure to approximate the distribution of this residual process, and she showed that the smooth bootstrap `works' in the sense that the limiting distribution of the bootstrapped residual process, conditional on the data, equals the process $W(y)$ defined above, in probability.  We will study an alternative bootstrap procedure that has the advantage of not requiring smoothing of the residual distribution.   For $i=1,\ldots,n$, let $\tilde \eps_i = \hat\eps_i - n^{-1} \sumj \hat\eps_j$, and let 
$$ \hat F_{0,n}(y) =  n^{-1} \sumi I\{\tilde\eps_i \le y\} $$
be the (non-smoothed) empirical distribution of the centered residuals.  Then, we randomly draw bootstrap errors $\eps_{0,1}^*,\ldots,\eps_{0,n}^*$ with replacement from $\hat F_{0,n}$.   Let $Y_i^* = \hat m(X_i) + \eps_{0,i}^*$, $i=1,\ldots,n$, and let $\hat m_0^*(\cdot)$ be the same as $\hat m(\cdot)$, except that we use the bootstrap data $(X_1,Y_1^*),\ldots,(X_n,Y_n^*)$. Define now
\begin{equation}\label{hateps-boot} 
\hat\eps_{0,i}^* = Y_i^* - \hat m_0^*(X_i) = \eps_{0,i}^* + \hat m(X_i) - \hat m_0^*(X_i). 
\end{equation}
We are interested in the asymptotic behavior of the process $n^{1/2} (\hat F_{0,n}^* - \hat F_{0,n})$ with
\begin{equation}\label{hat-F_0*} 
 \hat F_{0,n}^*(y) := n^{-1} \sumi I\{\hat\eps_{0,i}^* \le y\} 
\end{equation}
and we will show below that it converges asymptotically to the same limiting Gaussian process as the original residual process $n^{-1/2} \sumi \big(I\{\hat\eps_i \le y\} - F(y)\big)$, $y \in \R$, which means that smoothing of the residuals is not necessary to obtain a consistent bootstrap procedure.

In order to prove this result, we will use the results proved in Neumeyer (2009) to show that the difference between the smooth and the non-smooth bootstrap residual process is asymptotically negligible.   To this end, first note that we can write $\eps_{0,i}^* = \hat F_{0,n}^{-1}(U_i)$, $i=1,\ldots,n$, where $U_1, \ldots, U_n$ are i.i.d.\ random variables from a $U[0,1]$ distribution.  Strictly speaking the $U_i$'s form a triangular array $U_{1,n},\ldots,U_{n,n}$ of $U[0,1]$ variables, but since we are only interested in convergence in distribution of the bootstrap residual process (as opposed to convergence in probability or almost surely), we can work with $U_1,\ldots,U_n$ without loss of generality.  

We introduce the following notations : let $\eps_{s,i}^* = \hat F_{s,n}^{-1}(U_i)$, where $\hat F_{s,n}(y) = \int \hat F_{0,n}(y-vs_n) \, dL(v)$ is the convolution of the distribution $\hat F_{0,n}(y-\cdot s_n)$ and the integrated kernel $L(\cdot) = \int_{-\infty}^\cdot \ell(u) \, du$, where $\ell$ is a kernel density function and $s_n$ is a bandwidth sequence controlling the smoothness of $\hat F_{s,n}$.   Then, similarly to the definition of $\hat\eps_{0,i}^*$ in (\ref{hateps-boot}), we define 
\begin{equation}\label{hateps-boot-s0} 
\hat \eps_{s,0,i}^* = \eps_{s,i}^* + \hat m(X_i) - \hat m_0^*(X_i).
\end{equation}
 We then decompose the bootstrap residual process as follows :
\begin{eqnarray}
n^{1/2} \big(\hat F_{0,n}^*(y) - \hat F_{0,n}(y)\big) & = & n^{-1/2} \sumi \big(I\{\hat\eps_{0,i}^* \le y\} - I\{\hat\eps_{s,0,i}^* \le y\}\big) \nonumber \\
&& + n^{-1/2} \sumi \big(I\{\hat\eps_{s,0,i}^* \le y\} - \hat F_{s,n}(y) \big) \nonumber \\
&& + n^{1/2} \big(\hat F_{s,n}(y) - \hat F_{0,n}(y)\big) \nonumber \\[.1cm]
& = & T_{n1}(y) + T_{n2}(y) + T_{n3}(y). \label{decomp}
\end{eqnarray}
In Lemmas \ref{lem1} and \ref{lem3} in the Appendix we show that under assumptions \ref{A.1}--\ref{A.4} and conditions \ref{C.1}, \ref{C.2} concerning the choice of $\ell$ and $s_n$ in the proof (also given in the Appendix),  the terms $T_{n1}$ and $T_{n3}$ are asymptotically negligible. For the proof of negligibility of  $T_{n1}$ it is important that in (\ref{hateps-boot-s0}) the same function $\hat m_0^*$ is used than in (\ref{hateps-boot}) (in contrast to (\ref{hateps-boot-s}) below). Further in  Lemma \ref{lem2} we show that the process $T_{n2}$ is asymptotically equivalent to the smooth bootstrap residual process $n^{1/2} (\hat F_{s,n}^* - \hat F_{s,n})$ with 
\begin{eqnarray}\label{hat-F_s*}
\hat F_{s,n}^*(y) & = & n^{-1} \sumi I\{\hat\eps_{s,i}^* \le y\} ,
\end{eqnarray}
where, in contrast to (\ref{hateps-boot-s0}),
\begin{equation}\label{hateps-boot-s} 
\hat \eps_{s,i}^* = \eps_{s,i}^* + \hat m(X_i) - \hat m_s^*(X_i)
\end{equation}
with $\hat m_s^*$ defined as $\hat m$, but based on smoothed bootstrap data $(X_i,\hat m(X_i)+\eps_{s,i}^*)$, $i=1,\dots,n$. 
Neumeyer (2009) showed weak convergence of the residual process based on the smooth residual bootstrap, $n^{1/2} (\hat F_{s,n}^* - \hat F_{s,n})$, to the Gaussian process  defined in (\ref{covW}).     This leads to the following main result regarding the validity of the non-smooth bootstrap residual process.  

\begin{Theorem} \label{mainth}
Assume \ref{A.1}--\ref{A.4} in the Appendix.  Then, conditionally on the data $(X_1,Y_1),\ldots,(X_n,Y_n)$, the process 
$n^{1/2} \big(\hat F_{0,n}^*(y) - \hat F_{0,n}(y)\big)$, $y \in \R$, converges weakly to the zero-mean Gaussian process $W(y)$, $y \in \R$, defined in (\ref{covW}), in probability. 
\end{Theorem}

The proof is given in the Appendix.

\section{Linear model}  \label{sec-lin}

Consider independent observations from the linear model
\begin{eqnarray}
Y_{ni} &=& \bx_{ni}^T\bbeta+\eps_{ni},\quad i=1,\ldots,n, \label{mod}
\end{eqnarray}
where $\bbeta\in\R^p$ denotes the unknown parameter and the errors $\eps_{ni}$ are
assumed to be independent and identically distributed
with $E[\eps_{ni}]=0$ and distribution function $F$. Throughout this section let 
$\bX_n \in \R^{n\times p}$ denote the design matrix in the linear model, 
where the vector $\bx_{ni}^T=(x_{ni1},\ldots,x_{nip})$ corresponds to the $i$th row of the matrix $\bX_n$ and is
not random. The design matrix $\bX_n\in\R^{n\times p}$ is assumed to be of rank $p\leq n$ and to satisfy the  
regularity assumptions \ref{AL.1} given in the Appendix. 
We consider the least squares estimator 
\begin{eqnarray}\label{M}
\hat{\bbeta}_n\;=\;(\bX_n^T\bX_n)^{-1}\bX_n^T\bY_n\;=\; \bbeta+(\bX_n^T\bX_n)^{-1}\bX_n^T\beps_n
\end{eqnarray}
with the notations $\bY_n=(Y_{n1},\ldots,Y_{nn})^T,\beps_n=(\varepsilon_{n1},\ldots,\varepsilon_{nn})^T$, and define residuals $\hat\eps_{ni}=Y_{ni}-\bx_{ni}^T\hat\bbeta_n$, $i=1,\dots,n$.


Residual processes in linear models have been extensively studied by Koul (2002). It is shown there that, under our assumptions \ref{AL.1} and \ref{AL.2} in the Appendix, the process
$n^{-1/2} \sumi \big(I\{\hat\eps_{ni} \le y\} - F(y)\big)$, $y \in \R$, converges weakly to a zero-mean Gaussian process $ W(y)$ with covariance function 
\begin{eqnarray}\nonumber 
 Cov\big( W(y_1), W(y_2)\big) &=& F(y_1\wedge y_2)-F(y_1)F(y_2) +\bold{m}^T\boldsymbol{\Sigma}^{-1}\bold{m}\Big( f(y_1)f(y_2)Var(\eps)\\
&&{}+f(y_1)E[I\{\eps\leq y_2\}\eps]+f(y_2)E[I\{\eps\leq y_1\}\eps]\Big), 
\label{covW-lin}
\end{eqnarray}
where $\eps$ has distribution function $F$ and density $f$, and $\bold{m}$ and $\boldsymbol{\Sigma}$ are defined in \ref{AL.1}.

For the bootstrap procedure we generate $\eps_{0,i}^*$, $i=1,\dots,n$, from the distribution function 
$$ \hat F_{0,n}(y) =  n^{-1} \sumi I\{\hat\eps_{ni} \le y\}. $$
Note that unlike in the nonparametric case, we don't have to center the residuals, as they are centered by construction.  
Note also that in the bootstrap residuals we suppress the index $n$ to match the notation in the nonparametric case. 
We now define bootstrap observations by
\begin{eqnarray*}\label{boot}
Y_{ni}^*&=&\bx_{ni}^T\hat{\bbeta}_n+\eps_{0,i}^* \quad (i=1,\ldots,n),
\end{eqnarray*}
and calculate estimated residuals from the bootstrap sample 
\begin{equation}\label{hateps-boot-lin}
\hat{\eps}_{0,i}^* = Y_{ni}^*-\bx_{ni}^T\hat{\bbeta}_{0,n}^*=\eps_{0,i}^*+\bx_{ni}^T(\hat\bbeta_n-\hat{\bbeta}_{0,n}^*),
\end{equation}
where $\hat{\bbeta}_{0,n}^*$ is the least squares estimator
\begin{eqnarray}\label{hat-beta*}
\hat{\bbeta}_{0,n}^*\;=\;(\bX_n^T\bX_n)^{-1}\bX_n^T\bY_n^*\;=\; \hat\bbeta_n+(\bX_n^T\bX_n)^{-1}\bX_n^T\beps_{0,n}^*
\end{eqnarray}
(with the notations $\bY_n^*=(Y_{n1}^*,\ldots,Y_{nn}^*)^T,\beps_{0,n}^*=(\varepsilon_{0,1}^*,\ldots,\varepsilon_{0,n}^*)^T$).
We will show that the bootstrap residual process $n^{1/2} \big(\hat F_{0,n}^*(y) - \hat F_{0,n}(y)\big)$, with 
\begin{eqnarray} \label{hatF0nstar}
\hat F_{0,n}^*(y) = n^{-1} \sumi I\{\hat\eps_{0,i}^* \le y\}, 
\end{eqnarray}
converges to the same limiting process $W(y)$, $y\in\R$, as the original residual process $n^{-1/2} \sumi \big(I\{\hat\eps_{ni} \le y\} - F(y)\big)$, $y \in \R$.
Using the representations $\eps_{0,i}^* = \hat F_{0,n}^{-1}(U_i)$, $\eps_{s,i}^* = \hat F_{s,n}^{-1}(U_i)$, $i=1,\ldots,n$, where $U_1, \ldots, U_n$ are i.i.d.\ $U[0,1]$-distributed and   $\hat F_{s,n}(y) = \int \hat F_{0,n}(y-vs_n) \, dL(v)$ is the smoothed empirical distribution function of the residuals,  we have the same decomposition (\ref{decomp}) as in the nonparametric case. 
In Lemmas \ref{lem1-linear} and \ref{lem3-linear} in the Appendix we show that under assumptions \ref{AL.1}, \ref{AL.2} and conditions \ref{CL.1}, \ref{CL.2} on the choice of $s_n$ and $L$ (also given in the Appendix) the terms $T_{n1}$ and $T_{n3}$ are asymptotically negligible. Further we show in Lemma \ref{lem2-linear} that the limiting distribution of 
\begin{equation}\label{T_n2-lin}
T_{n2}(y)=n^{-1/2} \sumi \big(I\{\hat\eps_{s,0,i}^* \le y\} - \hat F_{s,n}(y) \big)
\end{equation}
with $\hat\eps_{s,0,i}^*=\eps_{s,i}^*+\bx_{ni}^T(\hat\bbeta_n-\hat\bbeta_{0,n}^*)$ (analogous to (\ref{hateps-boot-lin}), but with smooth bootstrap errors $\eps_{s,i}^*$)  
is the same as the limiting distribution of $n^{1/2} \big(\hat F_{s,n}^*(y) - \hat F_{s,n}(y) \big)$, with
\begin{eqnarray} \label{hatFsnstar}
\hat F_{s,n}^*(y) = n^{-1} \sumi I\{\hat\eps_{s,i}^* \le y\}
\end{eqnarray} 
and $\hat\eps_{s,i}^*=\eps_{s,i}^*+\bx_{ni}^T(\hat\bbeta_n-\hat\bbeta_{s,n}^*)$, where $\hat\bbeta_{s,n}^*=\hat\bbeta_n+(\bX_n^T\bX_n)^{-1}\bX_n^T\beps_{s,n}^*$ and $\beps_{s,n}^*=(\varepsilon_{s,1}^*,\ldots,\varepsilon_{s,n}^*)^T$.  
To show this we use results from Koul and Lahiri (1994). In this way we obtain the validity of the classical residual bootstrap.   

\begin{Theorem} \label{mainth-linear}
Assume \ref{AL.1}, \ref{AL.2}  in the Appendix.  Then, conditionally on the data $Y_{1n},\dots,Y_{nn}$, the process $n^{1/2} \big(\hat F_{0,n}^*(y) - \hat F_{0,n}(y)\big)$, $y \in \R$, converges weakly to the zero-mean Gaussian process $ W(y)$, $y \in \R$, defined in (\ref{covW-lin}), in probability. 
\end{Theorem}

The proof is given in the Appendix.

\section{Simulations}  \label{simulations}

In this section we will study the behavior of the smooth and the non-smooth residual bootstrap for a range of models, sample sizes and contexts.  We start with an empirical study to assess the quality of the approximation of the distribution of the residual empirical distribution by means of the bootstrap.  Next, we will consider the use of the residual bootstrap for approximating the critical values in two representative examples of hypothesis testing procedures based on residual processes : the case of testing for symmetry of the error distribution, and the case of goodness-of-fit tests for the regression function. 

\subsection{Bootstrap approximation}

Consider model (\ref{model}) in the nonparametric case and generate data with $m(x) = 2x$, where $X$ follows a uniform distribution on $[0,1]$, and $\eps \sim N(0,0.25^2)$.  In order to assess the quality of the smooth and the non-smooth bootstrap approximation, we compare the distribution of the least squares distance (or $L^2$ or Cr\'amer-von Mises distance)  $LS = \sum_{i=1}^n \big(\hat F_{0,n}(\hat\eps_i)-F(\hat\eps_i)\big)^2$ with its respective bootstrap versions based on the smooth or non-smooth residual bootstrap (see (\ref{hat-F_0*}) and (\ref{hat-F_s*}), respectively):
$$ LS_0^* = \sum_{i=1}^n \big(\hat F_{0,n}^*(\hat\eps_{0,i}^*) - \hat F_{0,n}(\hat\eps_{0,i}^*)\big)^2, \hspace*{.5cm} LS_s^* = \sum_{i=1}^n \big(\hat F_{s,n}^*(\hat\eps_{s,i}^*) - \hat F_{s,n}(\hat\eps_{s,i}^*)\big)^2. $$
Here, $\hat F_{0,n}$, $\hat F_{s,n}$,  $\hat F_{0,n}^*$ and $\hat F_{s,n}^*$ are defined as in Section \ref{nonpreg}. 
We also calculate the more robust median absolute deviation distance defined by $MAD = n \mbox{ median} \big(\big|\hat F_{0,n}(\hat\eps_i)-F(\hat\eps_i)\big|\big)_{i=1}^n$ and compare this with the bootstrap versions :
$$ MAD_0^* = n \mbox{ median} \big(\big|\hat F_{0,n}^*(\hat\eps_{0,i}^*) - \hat F_{0,n}(\hat\eps_{0,i}^*)\big|\big)_{i=1}^n, $$
$$ MAD_s^* = n \mbox{ median} \big(\big|\hat F_{s,n}^*(\hat\eps_{s,i}^*) - \hat F_{s,n}(\hat\eps_{s,i}^*)\big|\big)_{i=1}^n. $$
In order to verify whether the distribution of $LS$ is close to that of its bootstrap versions $LS_0^*$ and $LS_s^*$ (and similarly for $MAD$), we calculate the proportion of samples for which $LS$ exceeds the $(1-\alpha)$-th quantile of the distribution of $LS_0^*$ and $LS_s^*$.  If the bootstrap approximation works well, we expect a proportion close to $\alpha$.  Table \ref{table1.1} shows the results of this comparison for several sample sizes and several quantile levels $\alpha$.   The results are based on 500 simulation runs, and for each simulation 500 bootstrap samples are generated.   The bandwidth $h_n$ is taken equal to $h_n=\hat\sigma_X n^{-0.3}$ with $\hat\sigma_X$ the empirical standard deviation of $X_1,\ldots,X_n$, and the kernel $k$ is the biweight kernel.  For the smooth bootstrap, bootstrap errors $\eps_{s,i}^*$ are generated from $\hat F_{s,n}(y) = \int \hat F_{0,n}(y-vs_n) \, dL(v)$, where $L$ is the distribution of a standard normal and $s_n = 0.5 n^{-1/4}$, as in Neumeyer (2009).  

The table shows that the smooth bootstrap provides a quite good approximation especially for larger sample sizes, whereas the proportions provided by the non-smooth bootstrap are systematically too small.  The MAD distance provides a better approximation than the less robust LS distance.  Overall, it is clear that the non-smooth residual bootstrap, although it is asymptotically valid, performs less good than the smooth residual bootstrap in this situation.   This was already shown by Neumeyer (2009), and so here we confirm these findings using different models and different statistics, and taking a larger range of sample sizes.  A possible explanation for the discrepancy between the behavior of the smooth and the non-smooth residual bootstrap is that $F(y)$ (in the original statistic) and $\hat F_{s,n}(y)$ (in the smooth bootstrap approximation) are both smooth, whereas $\hat F_{0,n}(y)$ (in the non-smooth bootstrap) is a step function, which leads to a less good approximation.  

\begin{table}[!h]
\begin{center}
\hspace*{-.2cm}
\begin{tabular}{c|c|ccccccc}
$n$ & Distance & \multicolumn{7}{c}{$\alpha$}  \\
& & .025 & .05 & .1 & .3 & .5 & .7 & .9 \\
\hline
50    & $LS_s^*$     & .012 & .026 & .054 & .178 & .368 & .634 & .864 \\
        & $LS_0^*$     & .008 & .018 & .020 & .086 & .176 & .310 & .608 \\
        & $MAD_s^*$ & .010 & .024 & .044 & .156 & .346 & .590 & .856 \\
        & $MAD_0^*$ & .012 & .012 & .062 & .112 & .158 & .604 & .752 \\ \hline
100  & $LS_s^*$     & .012 & .030 & .072 & .258 & .434 & .634 & .874 \\
        & $LS_0^*$     & .004 & .004 & .028 & .124 & .272 & .398 & .652 \\
        & $MAD_s^*$ & .006 & .026 & .060 & .226 & .384 & .618 & .870 \\
        & $MAD_0^*$ & .006 & .018 & .058 & .152 & .350 & .358 & .744 \\ \hline
200  & $LS_s^*$     & .026 & .042 & .080 & .244 & .434 & .652 & .898 \\
        & $LS_0^*$     & .012 & .028 & .046 & .156 & .260 & .424 & .692 \\
        & $MAD_s^*$ & .022 & .044 & .078 & .232 & .432 & .662 & .906 \\
        & $MAD_0^*$ & .026 & .042 & .084 & .216 & .300 & .570 & .766 \\ \hline
500  & $LS_s^*$     & .014 & .036 & .072 & .246 & .458 & .680 & .880 \\
        & $LS_0^*$     & .008 & .024 & .050 & .170 & .324 & .496 & .734 \\
        & $MAD_s^*$ & .014 & .034 & .070 & .246 & .442 & .638 & .880 \\
        & $MAD_0^*$ & .014 & .038 & .062 & .244 & .422 & .464 & .782 \\ \hline
1000 & $LS_s^*$    & .020 & .044 & .090 & .288 & .482 & .692 & .916 \\
        & $LS_0^*$     & .018 & .034 & .064 & .226 & .390 & .544 & .786 \\
        & $MAD_s^*$ & .020 & .034 & .108 & .308 & .490 & .684 & .902 \\
        & $MAD_0^*$ & .020 & .038 & .104 & .292 & .422 & .598 & .856 \\ 
\end{tabular}
\caption{{\it Comparison between the behavior of the smooth and the non-smooth bootstrap in the nonparametric model for several sample sizes $n$, quantile levels $\alpha$ and distance measures.}}
\label{table1.1}
\vspace*{-.5cm}
\end{center}
\end{table}

\subsection{Application to testing for symmetry}

As already mentioned before, the residual bootstrap is very much used in hypothesis testing regarding various aspects of model (\ref{model}).  As a first illustration we consider a test for the symmetry of the error density in a linear regression model with fixed design.  More precisely, consider the model $Y_{ni} = \bx_{ni}^T \bbeta + \eps_{ni}$, where $E(\eps_{ni})=0$, and suppose we are interested in testing the following hypothesis regarding the distribution $F$ of $\eps_{ni}$ :
$$ H_0 : F(t) = 1 - F(-t) \mbox{ for all } t \in \R. $$
When the design is fixed and the regression function is linear, Koul (2002) considered a test for $H_0$ based on the residual process 
$$ \hat F_{0,n}(\cdot) - \hat F_{-0,n}(\cdot) = n^{-1} \sum_{i=1}^n \big[I(\hat\eps_{ni} \le \cdot) - I(-\hat\eps_{ni} \le \cdot)\big], $$ 
where $\hat F_{-0,n}$ is the empirical distribution of $-\hat\eps_{n1},\ldots,-\hat\eps_{nn}$, and $\hat\eps_{ni} = Y_{ni}-\bx_{ni}^T\hat{\bbeta}_n$.  Natural test statistics are the Kolmogorov-Smirnov and the Cr\'amer-von Mises statistics :
$$ T_{KS} = \sup_y |\hat F_{0,n}(y) - \hat F_{-0,n}(y)|, \:\:\: T_{CM} = \int \big(\hat F_{0,n}(y) - \hat F_{-0,n}(y)\big)^2 d\hat F_{0,n}(y). $$
It is clear from the covariance function given in (\ref{covW-lin}) that their asymptotic distribution is not easy to approximate, and that the residual bootstrap offers a valid alternative.   We will compare the level and power of the two tests, using the smooth and the non-smooth bootstrap.   The bootstrapped versions of $T_{CM}$ (and similarly for $T_{KS}$) are given by 
$$ T_{CM,0}^* = n \int \big(\hat F_{0,n}^*(y) - \hat F_{-0,n}^*(y)\big)^2 d\hat F_{0,n}^*(y), \:\:\: T_{CM,s}^* =n \int \big(\hat F_{s,n}^*(y) - \hat F_{-s,n}^*(y)\big)^2 d\hat F_{s,n}^*(y), $$ 
where for the non-smooth bootstrap, bootstrap errors $\eps_{0,1}^*,\ldots,\eps_{0,n}^*$ are drawn from \linebreak $\frac12 \big(\hat F_{0,n}(\cdot) + \hat F_{-0,n}(\cdot)\big)$, which is by construction a symmetric distribution, and for the smooth bootstrap we smooth this distribution using $s_n=0.5 n^{-1/4}$ and a Gaussian kernel.  The estimators $\hat F_{0,n}^*(\cdot)$ and $\hat F_{s,n}^*(\cdot)$ are defined as in (\ref{hatF0nstar}) and (\ref{hatFsnstar}), and $\hat F_{-0,n}^*(\cdot)$ and $\hat F_{-s,n}^*(\cdot)$ are defined accordingly.  Finally, we reject $H_0$ if the observed value of $T_{CM}$ exceeds the quantile of level $1-\alpha$ of the distribution of $T_{CM,0}^*$ or $T_{CM,s}^*$.

\begin{table}[!h]
\begin{center}
\hspace*{-.2cm}
\begin{tabular}{c|c|ccc|ccc|ccc}
$n$ & Test & \multicolumn{3}{c|}{$d=0$} & \multicolumn{3}{c|}{$d=2$} & \multicolumn{3}{c}{$d=4$} \\
& & .025 & .05 & .1 & .025 & .05 & .1 & .025 & .05 & .1 \\
\hline
50    & $T_{KS,s}^*$   & .019 & .035 & .060 & .051 & .087 & .139 & .139 & .207 & .284 \\
	& $T_{KS,0}^*$  & .012 & .026 & .054 & .038 & .066 & .128 & .116 & .175 & .266 \\
	& $T_{CM,s}^*$ & .027 & .049 & .103 & .081 & .132 & .210 & .202 & .292 & .413 \\
	& $T_{CM,0}^*$ & .021 & .044 & .093 & .067 & .117 & .186 & .178 & .267 & .389 \\ \hline
100  & $T_{KS,s}^*$  & .017 & .036 & .077 & .086 & .135 & .219 & .286 & .383 & .513 \\
	& $T_{KS,0}^*$  & .013 & .029 & .065 & .071 & .118 & .200 & .264 & .361 & .490 \\
	& $T_{CM,s}^*$ & .024 & .047 & .097 & .125 & .194 & .283 & .377 & .495 & .643 \\
	& $T_{CM,0}^*$ & .021 & .041 & .093 & .114 & .180 & .272 & .363 & .472 & .623 \\ \hline
200  & $T_{KS,s}^*$  & .024 & .044 & .088 & .179 & .269 & .384 & .598 & .710 & .823 \\
	& $T_{KS,0}^*$  & .019 & .036 & .072 & .160 & .241 & .357 & .570 & .688 & .804 \\
	& $T_{CM,s}^*$ & .026 & .047 & .097 & .228 & .334 & .462 & .727 & .826 & .908 \\
	& $T_{CM,0}^*$ & .024 & .047 & .092 & .222 & .319 & .447 & .718 & .820 & .902 \\ \hline
500  & $T_{KS,s}^*$  & .023 & .047 & .092 & .477 & .610 & .726 & .970 & .988 & .998 \\
	& $T_{KS,0}^*$  & .021 & .039 & .085 & .457 & .584 & .709 & .967 & .988 & .996 \\
	& $T_{CM,s}^*$ & .026 & .048 & .102 & .590 & .699 & .808 & .996 & .998 & .999 \\
	& $T_{CM,0}^*$ & .026 & .047 & .096 & .582 & .694 & .802 & .995 & .998 & .999 \\ \hline
1000 & $T_{KS,s}^*$ & .027 & .043 & .092 & .810 & .881 & .941 & 1	      & 1      &1 \\
	& $T_{KS,0}^*$ & .024 & .038 & .084 & .802 & .869 & .935 & 1	      & 1      & 1 \\
	& $T_{CM,s}^*$ & .025 & .049 & .097 & .897 & .946 & .975 & 1	      & 1      & 1 \\
	& $T_{CM,0}^*$ & .024 & .046 & .093 & .894 & .944 & .973 & 1	      & 1      & 1 \\ 
\end{tabular}
\vspace*{.3cm}
\caption{{\it Rejection probabilities of the test for symmetry in the linear model for several sample sizes $n$ and for $\alpha=0.025, 0.05$ and $0.1$.  Under the null we have a normal distribution ($d=0$), whereas under the alternative we have a skew normal distribution ($d=2$ and $d=4$). }}
\label{table3.1}
\vspace*{-.5cm}
\end{center}
\end{table}

Consider the model $Y_{ni} = 2x_{ni} + \eps_{ni}$, where $x_{ni} = i/n$.  We consider two error distributions under $H_0$.  The first one is a normal distribution with mean zero and variance $0.25^2$.  Under the alternative we consider the skew normal distribution of Azzalini (1985), whose density is given by $2\phi(y) \Phi(d y)$, where $\phi$ and $\Phi$ are the density and distribution of the standard normal.  More precisely, we let $d=2$ and $d=4$ and standardize these skew normal distributions so that they have mean zero and variance $0.25^2$.  Note that when $d=0$ we find back the normal distribution.    The second error distribution under $H_0$ is a Student-t distribution with 3 degrees of freedom, standardized in such a way that the variance equals $0.25^2$.  Note that the asymptotic theory does not cover this case, but we like to know how sensitive the bootstrap methods are to the existence of moments of higher order.   Under the alternative we consider a mixture of this Student-t distribution and a standard Gumbel distribution, again standardized to have mean zero and variance $0.25^2$.  The mixture proportions $p$ are 1 (corresponding to $H_0$), 0.75 and 0.50.  

The results, shown in Tables \ref{table3.1} and \ref{table3.2}, are based on 2000 simulation runs, and for each simulated sample a total of 2000 bootstrap samples are generated.   The tables show that the Cr\'amer-von Mises test outperforms the Kolmogorov-Smirnov test, and hence we focus on the former test.  Table \ref{table3.1} shows that for the normal error distribution, the level is about right for the smooth bootstrap and a little bit too low for the non-smooth bootstrap, resulting in a slightly higher power for the smooth bootstrap.  On the other hand, it should be noted that the behavior of the smooth bootstrap depends on the smoothing parameter $s_n$, whereas the non-smooth bootstrap is free of any additional parameters.   Table \ref{table3.2} shows a very different picture.   Here, the level is much too high for the smooth bootstrap, and about right for the non-smooth bootstrap.  When $n$ increases the situation for the smooth bootstrap gets better, but the level is still about 40-50$\%$ too high.  The power is somewhat higher for the smooth bootstrap, but given that the level is much too high, a fair comparison is not really possible here.   

\begin{table}[!h]
\begin{center}
\hspace*{-.2cm}
\begin{tabular}{c|c|ccc|ccc|ccc}
$n$ & Test & \multicolumn{3}{c|}{$p=1$} & \multicolumn{3}{c|}{$p=0.75$} & \multicolumn{3}{c}{$p=0.50$} \\
& & .025 & .05 & .1 & .025 & .05 & .1 & .025 & .05 & .1 \\
\hline
50    & $T_{KS,s}^*$  & .051 & .095 & .146 & .053 & .087 & .142 & .066 & .113 & .182 \\
	& $T_{KS,0}^*$ & .017 & .039 & .084 & .023 & .046 & .092 & .038 & .070 & .140 \\
	& $T_{CM,s}^*$ & .069 & .118 & .186 & .077 & .126 & .195 & .106 & .173 & .256 \\
	& $T_{CM,0}^*$ & .019 & .048 & .111 & .039 & .077 & .133 & .077 & .123 & .203 \\ \hline
100  & $T_{KS,s}^*$  & .058 & .093 & .165 & .067 & .103 & .164 & .113 & .170 & .266 \\
	& $T_{KS,0}^*$  & .021 & .044 & .085 & .033 & .060 & .113 & .074 & .121 & .206 \\
	& $T_{CM,s}^*$  & .069 & .112 & .184 & .088 & .131 & .211 & .167 & .243 & .338 \\
	& $T_{CM,0}^*$ & .028 & .053 & .106 & .046 & .078 & .139 & .109 & .182 & .283 \\ \hline
200  & $T_{KS,s}^*$  & .048 & .078 & .156 & .083 & .135 & .202 & .186 & .275 & .399 \\
	& $T_{KS,0}^*$  & .017 & .037 & .078 & .051 & .083 & .152 & .130 & .205 & .330 \\
	& $T_{CM,s}^*$ & .052 & .091 & .173 & .108 & .163 & .249 & .262 & .360 & .480 \\
	& $T_{CM,0}^*$ & .021 & .043 & .095 & .067 & .111 & .187 & .204 & .294 & .423 \\ \hline
500  & $T_{KS,s}^*$  & .041 & .075 & .135 & .131 & .200 & .288 & .438 & .553 & .672 \\
	& $T_{KS,0}^*$ & .019 & .038 & .086 & .082 & .143 & .235 & .383 & .485 & .619 \\
	& $T_{CM,s}^*$ & .046 & .086 & .151 & .169 & .248 & .349 & .558 & .653 & .756 \\
	& $T_{CM,0}^*$ & .021 & .043 & .100 & .116 & .192 & .296 & .496 & .607 & .716 \\ \hline
1000 & $T_{KS,s}^*$ & .041 & .070 & .136 & .226 & .330 & .434 & .776 & .860 & .909 \\
	& $T_{KS,0}^*$  & .020 & .044 & .089 & .165 & .257 & .375 & .728 & .818 & .895 \\
	& $T_{CM,s}^*$ & .045 & .079 & .136 & .290 & .385 & .498 & .857 & .907 & .949 \\
	& $T_{CM,0}^*$ & .023 & .048 & .100 & .230 & .328 & .446 & .837 & .890 & .939 \\ 
\end{tabular}
\vspace*{.3cm}
\caption{{\it Rejection probabilities of the test for symmetry in the linear model for several sample sizes $n$ and for $\alpha=0.025, 0.05$ and $0.1$.  Under the null we have a Student-t distribution with 3 degrees of freedom ($p=1$), whereas under the alternative we have a mixture of a Student-t (3) and a Gumbel distribution ($p=0.75$ and $p=0.50$).}}
\label{table3.2}
\vspace*{-.5cm}
\end{center}
\end{table}

\subsection{Application to goodness-of-fit tests}

We end this section with a second application of the residual bootstrap, which concerns the use of residual processes and the residual bootstrap for testing the fit of a parametric model for the regression function $m$ : 
$$ H_0 : m \in {\cal M} = \{m_\theta : \theta \in \Theta\}, $$
where ${\cal M}$ is a class of parametric regression functions depending on a $k$-dimensional parameter space $\Theta$.  Van Keilegom et al.\ (2008) showed that testing $H_0$ is equivalent to testing whether the error distribution satisfies $F \equiv F_0$, where $F_0$ is the distribution of $Y-m_{\theta_0}(X)$ and $\theta_0$ is the value of $\theta$ that minimizes $E[\{m(X)-m_\theta(X)\}^2]$. Consider the following Kolmogorov-Smirnov and Cr\'amer-von Mises type test statistics :
$$ T_{KS} = n^{1/2} \sup_y |\hat F_{0,n}(y)-\hat F_{\hat\theta}(y)|, \hspace*{1cm} T_{CM} = n \int \big(\hat F_{0,n}(y)-\hat F_{\hat\theta}(y)\big)^2 d\hat F_{\hat\theta}(y), $$
where $\hat F_{0,n}$ is as defined in Section \ref{nonpreg} and $\hat F_\theta(y) = n^{-1} \sum_{i=1}^n I(Y_i - \hat m_\theta(X_i) \le y)$ with 
$$ \hat m_{\theta}(x) = \sumi \frac{k_h(x-X_i)}{\sumj k_h(x-X_j)} m_{\theta}(X_i) $$
for any $\theta$, and $\hat\theta$ is the least squares estimator of $\theta$.  The critical values of these test statistics are approximated using our smooth and non-smooth residual bootstrap.  More precisely, the bootstrapped versions of $T_{CM}$ (and similarly for $T_{KS}$) are given by 
$$ T_{CM,0}^* = n \int \big(\hat F_{0,n}^*(y)-\hat F_{0,\hat\theta_0^*}^*(y)\big)^2 d\hat F_{0,\hat\theta_0^*}^*(y), $$
and 
$$ T_{CM,s}^* = n \int \big(\hat F_{s,n}^*(y)-\hat F_{s,\hat\theta_s^*}^*(y)\big)^2 d\hat F_{s,\hat\theta_s^*}^*(y),$$ 
where $\hat\theta_0^*$ is the least squares estimator based on the bootstrap data $(X_i,Y_{0,i}^*=m_{\hat\theta}(X_i)+\eps_{0,i}^*)$, $i=1,\ldots,n$, $\hat F_{0,\theta}^*(y) = n^{-1} \sum_{i=1}^n I(Y_{0,i}^* - \hat m_\theta(X_i) \le y)$ for any $\theta$, and similarly for $\hat\theta_s^*$ and $\hat F_{s,\theta}^*(y)$.  We reject $H_0$ if the observed value of $T_{CM}$ exceeds the quantile of level $1-\alpha$ of the distribution of $T_{CM,0}^*$ or $T_{CM,s}^*$.

\begin{table}[!h]
\begin{center}
\hspace*{-.2cm}
\begin{tabular}{c|c|ccc|ccc|ccc}
$n$ & Test & \multicolumn{3}{c|}{$a=0$} & \multicolumn{3}{c|}{$a=0.25$} & \multicolumn{3}{c}{$a=0.5$} \\
& & .025 & .05 & .1 & .025 & .05 & .1 & .025 & .05 & .1 \\
\hline
50    & $T_{KS,s}^*$  & .019 & .033 & .056 & .028 & .043 & .079 & .069 & .104 & .158 \\
        & $T_{KS,0}^*$  & .013 & .031 & .054 & .020 & .042 & .077 & .057 & .099 & .154 \\
        & $T_{CM,s}^*$ & .028 & .052 & .090 & .050 & .096 & .170 & .139 & .211 & .315 \\
        & $T_{CM,0}^*$ & .023 & .044 & .087 & .046 & .091 & .164 & .128 & .206 & .301 \\ \hline
100  & $T_{KS,s}^*$  & .015 & .033 & .065 & .052 & .077 & .129 & .142 & .200 & .303 \\
        & $T_{KS,0}^*$  & .011 & .025 & .054 & .035 & .068 & .116 & .115 & .180 & .279 \\
        & $T_{CM,s}^*$ & .026 & .054 & .102 & .089 & .135 & .227 & .274 & .376 & .512 \\
        & $T_{CM,0}^*$ & .021 & .047 & .090 & .078 & .121 & .211 & .249 & .356 & .494 \\ \hline
200  & $T_{KS,s}^*$  & .017 & .040 & .084 & .082 & .149 & .243 & .355 & .470 & .609 \\
        & $T_{KS,0}^*$  & .015 & .031 & .065 & .070 & .123 & .221 & .320 & .429 & .566 \\
        & $T_{CM,s}^*$ & .023 & .048 & .095 & .177 & .249 & .368 & .588 & .681 & .793 \\
        & $T_{CM,0}^*$ & .019 & .041 & .086 & .160 & .236 & .349 & .563 & .669 & .774 \\ \hline
500  & $T_{KS,s}^*$  & .015 & .032 & .071 & .247 & .350 & .454 & .819 & .894 & .942 \\
        & $T_{KS,0}^*$  & .011 & .023 & .063 & .213 & .318 & .433 & .795 & .876 & .939 \\
        & $T_{CM,s}^*$ & .022 & .040 & .091 & .395 & .509 & .630 & .947 & .970 & .983 \\
        & $T_{CM,0}^*$ & .019 & .038 & .083 & .379 & .495 & .613 & .941 & .968 & .982 \\ \hline
1000 & $T_{KS,s}^*$ & .020 & .042 & .088 & .517 & .641 & .748 & .996 & .998	& 1 \\
        & $T_{KS,0}^*$  & .018 & .037 & .079 & .488 & .609 & .734 & .994 & .998	& 1 \\
        & $T_{CM,s}^*$ & .023 & .047 & .099 & .700 &	.786 & .868 & .999 & 1 	& 1 \\
        & $T_{CM,0}^*$ & .020 & .044 & .092 & .685 & .774 & .864 & 1	     & 1	& 1 \\ 
\end{tabular}
\vspace*{.3cm}
\caption{{\it Rejection probabilities of the goodness-of-fit test for several sample sizes and for $\alpha=0.025, 0.05$ and $0.1$, when the error term has a normal distribution. The regression function is $m(x) = 2x + a x^2$ and the null hypothesis corresponds to $a=0$.}}
\label{table2.1}
\vspace*{-.5cm}
\end{center}
\end{table}

We consider the model $m(x) = 2x$ and let ${\cal M} = \{x \rightarrow \theta x : \theta \in \Theta\}$, i.e.\ the null model is a linear model without intercept.  The error term $\eps$ follows either a normal distribution or a Student-t distribution with 3 degrees of freedom, in both cases standardized in such a way that the variance equals $0.25^2$.  The covariate $X$ has a uniform distribution on $[0,1]$.  The bandwidth $h_n$ is taken equal to $h_n=\hat\sigma_X n^{-0.3}$, and the kernel $k$ is the biweight kernel.  For the smooth bootstrap, we use a standard normal distribution and $s_n = 0.5 n^{-1/4}$, as in Neumeyer (2009).   Under the alternative  we consider the model $m(x) = 2x + a x^2$ for $a=0.25$ and 0.5.  The rejection probabilities, given in Tables \ref{table2.1} and \ref{table2.2}, are based on 2000 simulation runs, and for each simulation 2000 bootstrap samples are generated.  

The tables show that the Cr\'amer-von Mises test outperforms the Kolmogorov-Smirnov test, independently of the sample size, the type of bootstrap and the value of $a$ (corresponding to the null or the alternative).  Hence, we focus on the Cr\'amer-von Mises test.  For the amount of smoothing chosen in this simulation (namely $s_n = 0.5 n^{-1/4}$) Table \ref{table2.1} shows that the smooth bootstrap behaves slightly better than the non-smooth bootstrap when the error is normal, but the difference is small.   However, when the error has a Student-t distribution (Table \ref{table2.2}) the level is too high for the smooth bootstrap, and more or less equal to the size of the test for the non-smooth bootstrap.  Similarly as in the case of tests for symmetry, a fair comparison of the power is not possible in this case, although the difference in power is nevertheless quite small. 

\begin{table}[!h]
\begin{center}
\hspace*{-.2cm}
\begin{tabular}{c|c|ccc|ccc|ccc}
$n$ & Test & \multicolumn{3}{c|}{$a=0$} & \multicolumn{3}{c|}{$a=0.25$} & \multicolumn{3}{c}{$a=0.5$} \\
& & .025 & .05 & .1 & .025 & .05 & .1 & .025 & .05 & .1 \\
\hline
50    & $T_{KS,s}^*$  & .020 & .034 & .063 & .056 & .080 & .123 & .134 & .182 & .256 \\
        & $T_{KS,0}^*$  & .016 & .030 & .057 & .044 & .069 & .117 & .117 & .167 & .250 \\
        & $T_{CM,s}^*$ & .026 & .048 & .098 & .102 & .152 & .229 & .243 & .334 & .450 \\
        & $T_{CM,0}^*$ & .021 & .044 & .089 & .091 & .142 & .223 & .231 & .320 & .436 \\ \hline
100  & $T_{KS,s}^*$  & .031 & .054 & .098 & .094 & .146 & .223 & .294 & .380 & .507 \\
        & $T_{KS,0}^*$  & .021 & .037 & .078 & .069 & .116 & .192 & .255 & .334 & .460 \\
        & $T_{CM,s}^*$ & .044 & .077 & .130 & .175 & .250 & .353 & .466 & .567 & .670 \\
        & $T_{CM,0}^*$ & .031 & .056 & .110 & .142 & .215 & .316 & .420 & .526 & .635 \\ \hline
200  & $T_{KS,s}^*$  & .038 & .064 & .109 & .179 & .248 & .360 & .583 & .695 & .797 \\
        & $T_{KS,0}^*$  & .023 & .041 & .082 & .135 & .197 & .300 & .509 & .617 & .746 \\
        & $T_{CM,s}^*$ & .044 & .078 & .139 & .283 & .363 & .474 & .751 & .825 & .888 \\
        & $T_{CM,0}^*$ & .026 & .053 & .103 & .226 & .315 & .427 & .703 & .781 & .859 \\ \hline
500  & $T_{KS,s}^*$  & .036 & .072 & .127 & .418 & .524 & .644 & .943 & .967 & .979 \\
        & $T_{KS,0}^*$  & .019 & .041 & .088 & .330 & .438 & .575 & .905 & .950 & .973 \\
        & $T_{CM,s}^*$ & .048 & .083 & .143 & .549 & .648 & .738 & .976 & .983 & .992 \\
        & $T_{CM,0}^*$ & .028 & .057 & .104 & .474 & .580 & .693 & .957 & .977 & .987 \\  \hline
1000 & $T_{KS,s}^*$ & .043 & .078 & .137 & .728 & .801 & .861 & .997 & .999 & .999\\
        & $T_{KS,0}^*$  & .028 & .049 & .095 & .638 & .744 & .824 & .994 & .998 & .999 \\
        & $T_{CM,s}^*$ & .047 & .079 & .132 & .796 & .853 & .910 & .998 & .999 & .999\\
        & $T_{CM,0}^*$ & .034 & .057 & .106 & .752 & .820 & .884 & .997 & .998 & .999\\ 
\end{tabular}
\vspace*{.3cm}
\caption{{\it Rejection probabilities of the goodness-of-fit test for several sample sizes and for $\alpha=0.025, 0.05$ and $0.1$, when the error term has a Student-t distribution with 3 degrees of freedom. The regression function is $m(x) = 2x + a x^2$ and the null hypothesis corresponds to $a=0$.}}
\label{table2.2}
\vspace*{-.5cm}
\end{center}
\end{table}

\appendix

\section{Proofs} \label{appendix}

\subsection{Nonparametric regression}

\noindent
The results of Section \ref{nonpreg} are valid under the following regularity conditions.

\begin{enumerate}[label=(\textbf{A\arabic{*}})]
\item\label{A.1} The univariate covariates $X_1,\dots,X_n$ are independent and identically distributed on a compact support, say $[0,1]$. They have a twice continuously differentiable density that is bounded away from zero. The regression function $m$ is twice continuously differentiable in $(0,1)$.
\item\label{A.2}  The errors $\eps_1,\dots,\eps_n$ are independent and identically distributed with distribution function $F$. They are centered and are independent of the covariates. $F$ is twice continuously differentiable with strictly positive density $f$ such that $\sup_{y\in\mathbb{R}}f(y)<\infty$ and $\sup_{y\in\mathbb{R}}|f^\prime(y)|<\infty$. Further, $E[|\eps_1|^{\upsilon}]<\infty$ for some $\upsilon\geq 7$.
\item\label{A.3}  $k$ is a twice continuously differentiable symmetric density with compact support $[-1,1]$, say, such that $\int uk(u)\,du=0$ and $k(-1)=k(1)=0$. The first derivative of $k$ is of bounded variation.
\item\label{A.4} $h_n$ is a sequence of positive bandwidths such that  $h_n\sim c_nn^{-\frac{1}{3}+\eta}$ with $\frac{4}{3+9\upsilon}<\eta<\frac{1}{12}$, $c_n$ is only of logarithmic rate, and $\upsilon$ is defined in \ref{A.2}. 
\end{enumerate}

Under those assumptions one in particular has that $nh_n^4=o(1)$ and it is possible to find some $\delta\in (0,\frac12)$ with 
\begin{eqnarray}\label{delta}
\frac{nh_n^{3+2\delta}}{\log(h_n^{-1})}\rightarrow \infty.
\end{eqnarray}
For the auxiliary smooth bootstrap residual process we need to choose a kernel $\ell$ and a bandwidth $s_n$ that fulfill the following conditions. 

\begin{enumerate}[label=(\textbf{C\arabic{*}})]
\item\label{C.1} Let $\ell$ denote a symmetric probability density function with compact support 
which has $\kappa\geq 2$ bounded derivatives inside the support, such that $\ell^{(\kappa)}\equiv 0$ and the first $\kappa-1$ derivatives are of bounded variation.
\item\label{C.2} Let $s_n=o(1)$ denote a positive sequence of bandwidths such that for $n\to\infty$,  
\begin{eqnarray*}
ns_n^4\to 0,\quad \frac{nh_ns_n^{2+8\delta/3}}{\log h_n^{-1}}\to\infty,\quad  h_n^{2\nu +2}=O(s_n^{2\nu-1}) \mbox{ for } \nu=1,\dots,\kappa-1, 
\end{eqnarray*}
where $\delta$ denotes the constant defined in (\ref{delta}). Further, $\delta>\frac{2}{\upsilon-1}$ with $\upsilon$ from \ref{A.2}. 
\end{enumerate}
 
As kernel $\ell$, e.g., the Epanechnikov kernel can be used. In order to fulfill \ref{C.2} one can choose $s_n\sim n^{-\frac14-\xi}$ for  $\xi>0$. 
Then $ns_n^4\to 0$ holds. Further with $h_n$ from assumption \ref{A.4} and (\ref{delta}),  the second bandwidth condition on $s_n$ in \ref{C.2} can be fulfilled (for very small $\xi$) if $\delta<\min(\frac32 (\eta+\frac16),\frac92\frac{\eta}{1-3\eta})$. This together with  $\delta >\frac{2}{\upsilon -1}$ gives the constraint on $\eta$ in assumption \ref{A.4}. The third bandwidth condition in \ref{C.2} is then fulfilled as well. 

Under assumptions \ref{A.1}--\ref{A.4} and conditions \ref{C.1}, \ref{C.2} assumptions (A.1)--(A.8) in Neumeyer (2009) (with the modification discussed in Remark 2) are satisfied (for smoothing parameter $a_n=s_n$ and kernel $k=\ell$).

Below we state and prove three lemmas regarding the terms $T_{n1}, T_{n2}$ and $T_{n3}$ given in (\ref{decomp}).  The proof of Theorem \ref{mainth} follows immediately from these three lemmas, together with Slutsky's Theorem.  

\begin{Lemma} \label{lem1}
Assume \ref{A.1}--\ref{A.4} and \ref{C.1}, \ref{C.2}.  Then, conditionally on the data, $\sup_y |T_{n1}(y)| = o_P(1)$, in probability, where $T_{n1}$ is defined in (\ref{decomp}).
\end{Lemma}

\noindent
{\bf Proof.} 
First note that by Lemma 3 in Cheng and Huang (2010), it is equivalent to show that $\sup_y |T_{n1}(y)|$ converges to zero in probability with respect to the joint probability measure of the original data and the bootstrap data.  Let $\hat g_n = (\hat m_0^*-\tilde m_n) - (\hat m-m_n)$ and $\kappa_n = \tilde m_n-m_n$, where 
$$ \tilde m_n(x) = \frac{\int K\big(\frac{t-x}{h_n}\big) m_n(t) f_X(t) \, dt}{\int K\big(\frac{t-x}{h_n}\big) f_X(t) \, dt}  \:\:\: \mbox{ and } \:\:\: m_n(x) = \frac{\int K\big(\frac{t-x}{h_n}\big) m(t) f_X(t) \, dt}{\int K\big(\frac{t-x}{h_n}\big) f_X(t) \, dt}, $$
and $f_X$ is the density of $X$.  Then,
\begin{eqnarray*}
T_{n1}(y) &=& n^{-1/2} \sumi \big(I\{\eps_{0,i}^* \le y + \hat g_n(X_i) + \kappa_n(X_i)\} - I\{\eps_{s,i}^* \le y + \hat g_n(X_i) + \kappa_n(X_i)\}\big) \\
&=& n^{-1/2} \sumi \big(I\{U_i \le \hat F_{0,n}(y + \hat g_n(X_i) + \kappa_n(X_i))\} - I\{U_i \le \hat F_{s,n}(y + \hat g_n(X_i) + \kappa_n(X_i))\}\big).
\end{eqnarray*}
Denoting $\hat w_n = \sup_y |\hat F_{s,n}(y) - \hat F_{0,n}(y)|$, it is clear that 
\begin{eqnarray} \label{upper}
&& T_{n1}(y) \le n^{-1/2} \sumi \big(I\{U_i \le \hat F_{s,n}(y + \hat g_n(X_i) + \kappa_n(X_i)) + \hat w_n\} \\
&& \hspace*{3.5cm} - I\{U_i \le \hat F_{s,n}(y + \hat g_n(X_i) + \kappa_n(X_i))\}\big), \nonumber
\end{eqnarray}
and similarly a lower bound for $T_{n1}(y)$ is given by the expression in which we replace $\hat w_n$ in the above upper bound by $-\hat w_n$.  Now consider the process 
$$ E_n(y,g,H,w) = n^{-1/2} \sumi \big(I\{U_i \le H(y + g(X_i) + \kappa_n(X_i)) + w\} - E[H(y+g(X) + \kappa_n(X))] - w\big), $$
indexed by $y \in \R, g \in {\cal G}, H \in {\cal H}$ and $w \in [-1,1]$.  Here, ${\cal G} = C_1^{1+\delta/2}([0,1])$ is the space of differentiable functions $g: [0,1] \rightarrow \R$ with derivatives $g'$ satisfying
\begin{eqnarray} \label{norm}
\max\Big\{\sup_{x\in [0,1]}|g(x)|,\sup_{x\in [0,1]}|g'(x)|\Big\}+\sup_{x_1,x_2 \in [0,1]}\frac{|g'(x_1)-g'(x_2)|}{|x_1-x_2|^{\delta/2}} & \le & 1,
\end{eqnarray}
see Van der Vaart and Wellner (1996, p.\ 154), with $\delta$ defined in (\ref{delta}). For $C=2\max\{1,\sup_y f(y) + \sup_y |f'(y)|\}$ let ${\cal H}$ denote the class of continuously differentiable distribution functions $H : \R \to [0,1]$ with uniformly bounded derivative $h$, such that 
\begin{eqnarray}\label{def-C}
\sup_{z \in \R} h(z) + \sup_{z_1,z_2 \in \R}\frac{|h(z_1)-h(z_2)|}{|z_1-z_2|^{\delta/2}} \le C,
\end{eqnarray}
and the tail condition 
\begin{eqnarray}\label{tailcond}
|1-H(z)| \le a/z^\upsilon \mbox{ for all $z \in \R^+$, } |H(z)|\le a/|z|^\upsilon \mbox{ for all $z \in \R^-$}
\end{eqnarray}
is satisfied for $\upsilon$ as in assumption \ref{A.2}, where the constant $a$ is independent of $H \in {\cal H}$.

In Proposition 3 in Neumeyer (2009) the asymptotic equicontinuity of the process $E_n(y,g,H,0)$ as a process in $(y,g,H) \in \R \times {\cal G} \times {\cal H}$ has been shown.  Using similar arguments the asymptotic equicontinuity of the extended process $E_n(y,g,H,w)$ in $\R \times {\cal G} \times {\cal H} \times [-1,1]$ can be shown, implying that $\sup_{y,g,H}|E_n(y,g,H,w) - E_n(y,g,H,0)| \rightarrow 0$ in probability as $w \rightarrow 0$.  Next, it follows from Lemma 2 and Proposition 4 in Neumeyer (2009)  that $P(\hat F_{s,n} \in {\cal H}) \rightarrow 1$, and we know from Proposition 2 and Lemma 3 in Neumeyer (2009) that $P(\hat g_n \in {\cal G}) \rightarrow 1$.  Note however that the estimator $\hat m_0^*$ in this paper is not exactly the same as the one in Neumeyer (2009), since our $\hat m_0^*$ is based on the non-smooth bootstrap errors $\eps_{0,i}^*$, whereas Neumeyer (2009) considers $\hat m_s^*$ as in (\ref{hateps-boot-s}),  based on the smooth bootstrap errors $\eps_{s,i}^*$, $i=1,\ldots,n$.  This means that in the proof of Lemma 3 in the latter paper the term depending on $s_n$ (or $a_n$ in the notation of that paper) does not exist.  So there is one term less to handle in our case, which means that the proof for $P(\hat g_n \in {\cal G}) \rightarrow 1$ in our case is actually simpler than in the latter paper.  Finally, $P(\hat w_n \in [-1,1]) \rightarrow 1$ by Lemma \ref{lem3} below.  This shows that the upper bound (\ref{upper}) is bounded by $n^{1/2}\hat w_n + o_P(1)$, uniformly in $y$, which is $o_P(1)$ by Lemma \ref{lem3} below.   Similarly, it can be shown that the lower bound is $o_P(1)$.  \hfill $\Box$

\begin{Lemma} \label{lem2}
Assume \ref{A.1}--\ref{A.4} and \ref{C.1}, \ref{C.2}.  Then, conditionally on the data, the process $T_{n2}(y)$, $y \in \R$, defined in (\ref{decomp}) converges weakly to the Gaussian process $W(y)$, $y \in \R$, in probability,  where $W$ was defined in (\ref{covW}).
\end{Lemma}

\noindent
{\bf Proof.} First note that the only difference between the process $T_{n2}$ and the process $R_n^*=n^{-1/2}(\hat F_{s,n}^*-\hat F_{s,n})$ (with $\hat F_{s,n}^*$ defined in (\ref{hat-F_s*})) studied in Neumeyer (2009) is the use of $\hat m_0^*$ instead of $\hat m_s^*$ (compare (\ref{hateps-boot-s0}) and (\ref{hateps-boot-s})). 
Hence, we should verify where the precise form of the estimator $\hat m_0^*$ is used in the proof of the weak convergence of $R_n^*$, which is given in Theorem 2 in Neumeyer's paper.   Carefully checking the proof of the latter theorem and of Lemma 1 in the same paper, reveals that $T_{n2}(y) = \tilde T_{n2}(y) + o_P(1)$, where 
$$ \tilde T_{n2}(y) = n^{-1/2} \sumi \big(I\{\eps_{s,i}^* \le y\} - \hat F_{s,n}(y) + \hat f_{s,n}(y) \eps_{0,i}^*\big), $$
where $\hat f_{s,n}(y) = (d/dy) \hat F_{s,n}(y)$, so $\tilde T_{n2}(y)$ equals $\tilde R_n^*(y)$ in the proof of Theorem 2 in Neumeyer (2009), except that one of the $\eps_{s,i}^*$'s has been replaced by $\eps_{0,i}^*$.  This has however no consequences for the proof of tightness and Lindeberg's condition in the proof of the aforementioned theorem. The only difference is that convergence of $E^*[I\{\eps_{s,i}^* \le y\}\eps_{0,i}^*]$ to $E[I\{\eps\le y\}\eps]$ in probability has to be shown. Here $E^*$ denotes conditional expectation, given the original data. 
Note that $E^*[I\{\eps_{s,i}^* \le y\}\eps_{0,i}^*]=E^*[I\{\eps_{s,i}^* \le y\}\eps_{s,i}^*]+r_n$, where $E^*[I\{\eps_{s,i}^* \le y\}\eps_{s,i}^*]$ converges to the desired $E[I\{\eps\le y\}\eps]$ according to Neumeyer (2009), and $|r_n|\leq E^*[|\eps_{0,i}^*-\eps_{s,i}^*|]$. For this we have 
\begin{eqnarray}\nonumber
E^*[|\eps_{0,i}^*-\eps_{s,i}^*|]&=&\int_0^1  |\hat F_{0,n}^{-1}(u)-\hat F_{s,n}^{-1}(u)|\,du\\
&\leq& \int_0^1  |\hat F_{0,n}^{-1}(u)-F^{-1}(u)|\,du +\int_0^1  |\hat F_{s,n}^{-1}(u)-F^{-1}(u)|\,du.
\label{wasserstein}
\end{eqnarray}
To conclude the proof we will show that the latter two integrals (which are the Wasserstein distance between $\hat F_{0,n}$ and $F$ and between $\hat F_{s,n}$ and $F$, respectively) converge to zero in probability. To this end let $\epsilon>0$. Note that $\int_0^1  |F^{-1}(u)|\,du=E[|\eps|]<\infty$ and let $\delta>0$ be so small that
\begin{eqnarray}\label{wasserstein1}
\int_{[\delta,1-\delta]^c}  |F^{-1}(u)|\,du &<& \epsilon.
\end{eqnarray}
Then note that 
\begin{eqnarray}\label{wasserstein2}
\int_{[\delta,1-\delta]}|\hat F_{0,n}^{-1}(u)-F^{-1}(u)|\,du &\leq& \sup_{u\in [\delta,1-\delta]} |\hat F_{0,n}^{-1}(u)-F^{-1}(u)| \;=\; o_P(1)
\end{eqnarray}
by the functional delta-method for quantile processes (see, e.g.\ Van der Vaart and Wellner (1996), p.\ 386/387) and uniform consistency of the residual empirical distribution function $\hat F_{0,n}$. Thus we have
\begin{eqnarray*}
 \int_0^1  |\hat F_{0,n}^{-1}(u)-F^{-1}(u)|\,du &\leq& \epsilon+\int_{[\delta,1-\delta]}|\hat F_{0,n}^{-1}(u)-F^{-1}(u)|\,du +\int_{[\delta,1-\delta]^c}  |\hat F_{0,n}^{-1}(u)|\,du\\
&\leq& \epsilon +o_P(1)+\Big| \int_0^1 (|\hat F_{0,n}^{-1}(u)|-| F^{-1}(u)|)\,du\Big| \\
&&{}+\Big| \int_{[\delta,1-\delta]} (|\hat F_{0,n}^{-1}(u)|-| F^{-1}(u)|)\,du\Big| + \int_{[\delta,1-\delta]^c}  |F^{-1}(u)|\,du\\
&\leq & 2\epsilon +o_P(1) +\Big| \int |x|\,d\hat F_{0,n}(x) -\int |x|\,d F(x)\Big|\\
&&{} +  \int_{[\delta,1-\delta]} |\hat F_{0,n}^{-1}(u)- F^{-1}(u)|\,du \\
&=& 2\epsilon +o_P(1),
\end{eqnarray*}
where the first inequality follows from (\ref{wasserstein1}), the second from (\ref{wasserstein2}), and the last equality follows from (\ref{wasserstein1}), (\ref{wasserstein2}) and convergence of $\int |x|\,d\hat F_{0,n}(x)=n^{-1}\sum_{i=1}^n |\hat\eps_i|$ to $E[|\eps|]=\int |x|\,d F(x)$ in probability. Analogous considerations for the second integral in (\ref{wasserstein}) make use of Lemma \ref{lem3} to obtain uniform consistency of $\hat F_{s,n}$, and of $\int |x|\,d\hat F_{0,n}(x)=n^{-1}\sum_{i=1}^n\int  |\hat\eps_i+s_nv|\ell(v)\,dv\to E[|\eps|]$ in probability. 
 \hfill $\Box$

\smallskip

\begin{Lemma} \label{lem3}
Assume \ref{A.1}--\ref{A.4} and \ref{C.1}, \ref{C.2}.  Then, conditionally on the data, $\sup_y |T_{n3}(y)| = o_P(1)$, in probability,  where $T_{n3}$ is defined in (\ref{decomp}).
\end{Lemma}

\noindent
{\bf Proof.} Write
\begin{eqnarray*}
\hat F_{s,n}(y) - \hat F_{0,n}(y) &=& \int \big[\hat F_{0,n}(y-vs_n) - \hat F_{0,n}(y)\big] \, dL(v) \\
& = & \int \big[\hat F_{0,n}(y-vs_n) - \hat F_{0,n}(y) - F(y-vs_n) + F(y) \big] \, dL(v) \\
&& + \int \big[F(y-vs_n) - F(y) \big] \, dL(v).
\end{eqnarray*}
The second term above is $O_P(s_n^2) = o_P(n^{-1/2})$ uniformly in $y$, since $f$ has a bounded derivative, $ns_n^4 \rightarrow 0$ and  $\int v \ell(v) \, dv = 0$.  For the first term we use the i.i.d.\ decomposition of $\hat F_{0,n}(y) - F(y)$ given in Theorem 1 in Akritas and Van Keilegom (2001), i.e.
\begin{eqnarray*}
\hat F_{0,n}(y) - F(y)& =& n^{-1} \sumi (I\{\eps_i\leq y\}  - F(y)) +f(y) n^{-1} \sumi \eps_i \\
\nonumber &&{}+ O(h_n^2) + o_P(n^{-1/2})
\end{eqnarray*}
uniformly in $y$.
Note that here we also apply $n^{-1} \sumi \hat\eps_i=o_P(n^{-1/2})$ which follows from M\"uller et al.\ (2004). 
 The expansion yields, uniformly in $y$,
\begin{eqnarray*}
&& \int \big[\hat F_{0,n}(y-vs_n) - \hat F_{0,n}(y) - F(y-vs_n) + F(y) \big] \, dL(v) \\
&& = n^{-1/2} \int \big[ E_n( y-vs_n) -E_n(y) \big] \, dL(v) + n^{-1} \sumi \eps_i \int \big[f(y-vs_n) - f(y)\big] \, dL(v)\\
&& \hspace*{.5cm}  {}+ O(h_n^2) + o_P(n^{-1/2})
\end{eqnarray*}
with the empirical process $E_n(y) = n^{-1/2} \sumi (I\{\eps_i \le y\} - F(y))$. 
From asymptotic equicontinuity of the process $E_n$ together with the bounded support of the kernel $\ell$ it follows that the first term is of order $o_P(n^{-1/2})$ uniformly in $y$. 
The second integral is of order $o_P(n^{-1/2})$ since $n^{-1} \sumi \eps_i = O_P(n^{-1/2})$, and $f$ is differentiable with bounded derivative. Further  $h_n^2=o(n^{-1/2})$ by the bandwidth conditions.
 \hfill $\Box$

\subsection{Linear model}

The results of Section \ref{sec-lin} are valid under the following regularity conditions.

\begin{enumerate}[label=(\textbf{AL\arabic{*}})]
\item\label{AL.1} The fixed design fulfills
\begin{enumerate}
\item\label{1} $\max_{i=1,\dots,n} \bx_{ni}^T(\bX_n^T\bX_n)^{-1}\bx_{ni}=O(\frac{1}{n})$,
\item\label{2} $\lim_{n\to\infty}\frac{1}{n}\bX_n^T\bX_n= \boldsymbol{\Sigma}\in\R^{p\times p}$ with invertible $\boldsymbol{\Sigma}$,
\item\label{3}$\lim_{n\to\infty}\frac{1}{n}\sum_{i=1}^n\bx_{ni}=\bold{m}\in\R^p$.
\end{enumerate}
\item\label{AL.2} The errors $\varepsilon_{ni}$, $i=1,\dots,n$, $n\in\mathbb{N}$, are independent and identically distributed with distribution function $F$ and density $f$ that is strictly positive, bounded, and continuously differentiable with bounded derivative on $\R$.  Assume $E[|\eps_1|^\upsilon]<\infty$ for some $\upsilon>3$. 
\end{enumerate}

The following conditions are needed for the auxiliary smooth bootstrap process.  
\begin{enumerate}[label=(\textbf{CL\arabic{*}})]
\item\label{CL.1} Let $\ell$ denote a symmetric and differentiable probability density function that is strictly positive on $\R$ with $\int u\ell(u)\,du=0$, $\int u^2\ell(u)\,du<\infty$. Let $\ell$ have $\kappa \ge 2$ bounded and square-integrable derivatives such that the first $\kappa-1$ derivatives are of bounded variation.
\item\label{CL.2} Let $s_n=o(1)$ denote a positive sequence of bandwidths such that 
$ns_n^4\to 0$ for $n\to\infty$.
For some $\delta\in (\frac{2}{\upsilon-1},2)$ with $\upsilon$ from \ref{AL.2} and for $\kappa$ from \ref{CL.1}, let $(\sqrt{n}s_n)^{\kappa(1-\delta/2)}s_n\to \infty$. 
\end{enumerate}

Note that condition \ref{CL.2} is possible if $\kappa$ is sufficiently large, namely $\kappa > 2/(2-\delta)$. 

\begin{Lemma} \label{lem1-linear}
Assume the linear model (\ref{mod}) and \ref{AL.1}, \ref{AL.2}, \ref{CL.1} and \ref{CL.2}.    Then, conditionally on the data, $\sup_y |T_{n1}(y)| = o_P(1)$, in probability, where
$$T_{n1}(y)= 
n^{-1/2} \sumi \big(I\{\hat\eps_{0,i}^* \le y\} - I\{\hat\eps_{s,0,i}^* \le y\}\big) . $$
\end{Lemma}

\noindent
{\bf Proof.} The proof is similar to the proof of Lemma \ref{lem1}. Note that 
\begin{eqnarray*}
T_{n1}(y) &=& n^{-1/2} \sumi \big(I\{U_i \le \hat F_{0,n}(y + \bx_{ni}^T(\hat\bbeta^*_{0,n}-\hat\bbeta_n)\} - I\{U_i \le \hat F_{s,n}(y + \bx_{ni}^T(\hat\bbeta^*_{0,n}-\hat\bbeta_n)\}\big).
\end{eqnarray*}
Thus with $\hat w_n = \sup_y |\hat F_{s,n}(y) - \hat F_{0,n}(y)|=o_P(n^{-1/2})$ (see Lemma \ref{lem3-linear}) one has the upper bound
\begin{eqnarray*}
&& T_{n1}(y) \le n^{-1/2} \sumi \big(I\{U_i \le \hat F_{s,n}(y + \bx_{ni}^T(\hat\bbeta^*_{0,n}-\hat\bbeta_n)) + \hat w_n\} \\
&& \hspace*{3.5cm} - I\{U_i \le \hat F_{s,n}(y + \bx_{ni}^T(\hat\bbeta^*_{0,n}-\hat\bbeta_n))\}\big), 
\end{eqnarray*}
and similarly a lower bound. 
Now write $\bx_{ni}^T(\hat\bbeta^*_{0,n}-\hat\bbeta_n)=\bz_{ni}^T{\bold b}_n$ 
with 
$$\bz_{ni}^T=n^{1/2}\bx_{ni}^T(\bX_n^T\bX_n)^{-1/2} \mbox{ and }\bb_n=n^{-1/2}(\bX_n^T\bX_n)^{1/2}(\hat\bbeta^*_{0,n}-\hat\bbeta_n).$$
Note that from \ref{1} the existence of a constant $K$ follows with $\| \bz_{ni}\|\leq K$ for all $i,n$. 
Further, from (\ref{hat-beta*}) we have $\bb_n=n^{-1/2}(\bX_n^T\bX_n)^{-1/2}\bX_n^T\beps_{0,n}^*$. Denoting by $Var^*$ the conditional variance, given the original data, it is easy to see that, for any $\eta>0$,
\begin{eqnarray*}
P(\|\bb_n\|>\eta) &\leq& \frac{1}{\eta^2} E\big[ Var^*( \bb_n )\big] =\frac{1}{n\eta^2} E\Big[\frac{1}{n}\sum_{j=1}^n \hat\eps_{nj}^2\Big] =o(1). 
\end{eqnarray*}
Thus we have $P(\bb_n\in B_\eta(0))\to 1$ for $n\to\infty$, where $B_\eta(0)$ denotes the ball of radius $\eta$ around the origin in $\R^p$. 

Let $\mathcal{H}$ denote the function class that was defined in the proof of Lemma \ref{lem1}. Note that the function class depends on $\upsilon$ from \ref{AL.2} and $\delta$ from \ref{CL.2}. In order to obtain $P(\hat F_{s,n} \in {\cal H}) \rightarrow 1$ one can mimic the proof of Lemma 2 (and Proposition 4) in Neumeyer (2009). In the linear case the proof actually is easier due to a faster rate of the regression estimator (one simply replaces $o(\beta_n)$ with $O(n^{-1/2})$ and sets $t_{ni}=0$). Then one obtains 
$\|\hat f_{n,s}-f\|_\infty=O(n^{-1/2}s_n^{-1})+O(n^{-\kappa/2}s_n^{\kappa+1})$ (where $\|\cdot\|_\infty$ is the supremum norm) and from this it follows that
$$\sup_{z_1\neq z_2}\frac{|\hat f_{n,s}(z_1)-f(z_1)-\hat f_{n,s}(z_2)+f(z_2)|}{|z_1-z_2|^{\delta/2}}=o(1)$$
under the bandwidth condition \ref{CL.2}. 

  Now consider the process $E_n=\sum_{i=1}^n (Z_{ni}-E[Z_{ni}])$ with 
$$ Z_{ni}(y,g,H,w) = n^{-1/2} I\{U_i \le H(y + g(\bz_{ni}) ) + w\} , $$
indexed by $y \in \R, g \in \tilde{\cal G}=\{g_\bb\mid \bb\in B_\eta(0)\}, H \in {\cal H}$ and $w \in [-1,1]$.  Here, $g_\bb: B_K(0)\to\R$, $g_\bb(\bz)=\bz^T\bb$, where $B_K(0)$ denotes the ball with radius $K$ around the origin in $\R^p$. One can apply Theorem 2.11.9 in Van der Vaart and Wellner (1996) in order to asymptotic equicontinuity of $E_n$. Here the main issue is to find a bound for the bracketing number $N_{[\,]}(\epsilon,{\cal F}, L_2^n)$, where ${\cal F}=\R\times\tilde {\cal G}\times{\cal H}\times[-1,1]$. This is the minimal number of sets in a partition of $\cal F$ into sets ${\cal F}_{\epsilon j}^n$ such that
\begin{equation}\label{brack}
\sum_{i=1}^n E\Big[\sup_{f_1,f_2\in{\cal F}_{\epsilon j}^n}|Z_{ni}(f_1)-Z_{ni}(f_2)|^2\Big]\leq \epsilon^2.
\end{equation}
To obtain a bound on this bracketing number, define $\tilde\epsilon=\epsilon/(3+C)^{1/2}$ with $C$ from (\ref{def-C}), and consider the following brackets for our index components.
\begin{itemize}
\item There are $O(\tilde\epsilon^{-2p})$ brackets of the form $[g_l,g_u]$ covering $\tilde{ \cal G}$ with  $\sup$-norm length $\tilde\epsilon^2$ according to Theorem 2.7.11 in Van der Vaart and Wellner (1996).
\item There are  $O(\exp(\tilde\epsilon^{-2/(1+\alpha)}))$ balls covering $\cal H$ with centers $H_c$ and  $\sup$-norm radius $\tilde\epsilon^2/2$ according to Lemma 4 in Neumeyer (2009).  Then one obtains as many brackets of the form $[H_l,H_u]=[H_c-\tilde\epsilon^2/2,H_c-\tilde\epsilon^2/2]$ that cover $\cal H$ and have sup-norm length $\tilde \epsilon^2$.  Here 
$\alpha= (\delta(\upsilon-1)/2-1)/(1+\upsilon+\delta/2) >0$ due to \ref{CL.2}.
\item Given $H_u$ and $g_u$ construct brackets $[y_l,y_u]$ such that $F_{n}(y_u|H_u,g_u)-F_{n}(y_l|H_u,g_u)\leq \tilde\epsilon^2$ for $F_{n}(y|H_u,g_u)=n^{-1}\sum_{i=1}^n H_u(y+g_u(\bz_{ni}))$. 
The brackets depend on $n$, but their minimal number is $O(\tilde\epsilon^{-2})$. 
\item There are $O(\tilde\epsilon^{-2})$ intervals $[w_l,w_u]$ of length $\tilde\epsilon^2$ that cover $[-1,1]$. 
\end{itemize}
The sets $[y_l,y_u]\times[g_l,g_u]\times[H_l,H_u]\times[w_l,w_u]$ (where $y_l,y_u$ correspond to $g_u$, $H_u$ as above) partition ${\cal F}$ and it is easy to see that (\ref{brack}) is fulfilled. We further obtain
$$N_{[\,]}(\epsilon,{\cal F}, L_2^n)=O(\epsilon^{-2p-4}\exp(\epsilon^{-2/(1+\alpha)})) $$
and thus the bracketing integral condition in Theorem 2.11.9 in Van der Vaart and Wellner (1996) holds. Further $\cal F$ is a totally bounded space with the semi-metric 
\begin{eqnarray*}
&&\rho((y_1,g_1,H_1,w_1),(y_2,g_2,H_2,w_2))\\
&=&\max\Big\{\sup_{g\in\tilde{\cal G}}\sup_{\bz\in B_K(0)}|H_1(y_1+g(\bz))-H_2(y_2+g(\bz))|,\|g_1-g_2\|_\infty, |w_1-w_2|\Big\}. 
\end{eqnarray*}
To see this we show that the covering number $N(\epsilon,{\cal F},\rho)$ is finite for every $\epsilon>0$. Let $(y,g,H,w)$ be some fixed arbitrary element of $\cal F$.
\begin{itemize}
\item There are finitely many balls with  $\sup$-norm radius $\epsilon$ covering $\tilde{ \cal G}$. Denote the centers of the balls as $g_1,\dots,g_{N_{\cal G}}$ and let $g_{c}^*$ denote the center of the ball containing $g$.   
\item There are  finitely many balls with  $\sup$-norm radius $\epsilon$ covering $\cal H$.  Denote the centers of the balls as $H_1,\dots,H_{N_{\cal H}}$ and let $H_{c}^*$ denote the center of the ball containing $H$. 
\item There are finitely many balls with radius $\epsilon$ covering $B_K(0)$ in $\R^p$. Denote the centers of the balls as $\bz_1,\dots,\bz_{N_{B}}$.  
\item Given $H_j$, $g_k$ and $\bz_l$, let the values $y_{j,k,l,m}$, $m=1,\dots,N_{j,k,l}$ segment the real line in finitely many intervals of length less than or equal to $\epsilon$ according to the distribution function  $y\mapsto H_j(y+g_k(\bz_l))$. Let $y_0=-\infty < y_1<\dots <y_{N_{R}}=\infty$ be the (finitely many) ordered values of all $y_{j,k,l,m}$, $j=1,\dots, N_{\cal H}$, $k=1,\dots, N_{\cal G}$, $l=1,\dots,N_{B}$. 

For $H_j=H_c^*$, $g_k=g_c^*$ let $y_c^*$ be the value $y_{j,k,l,m}$ closest to $y$, and, for all $l\in\{1,\dots,N_B\}$ denote the interval  $[y_{j,k,l,m-1},y_{j,k,l,m+1}]$ containing both $y$ and $y^*$ as $[y_{l,c}^{*1},y_{l,c}^{*2}]$. 
\item There are finitely many intervals $[ w_j-\epsilon/2,w_j+\epsilon/2]$, $j=1,\dots,N_{W}$ of length $\epsilon$ that cover $[-1,1]$. Let $w\in [ w_c-\epsilon/2,w_c+\epsilon/2]$.
\end{itemize}
Now we show that $(y,g,H,w)$ lies in the $(4+(1+\eta)C)\epsilon$-ball with respect to $\rho$ with center $(y_c^*,g_c^*,H_c^*,w_c^*)$. Applying the mean value theorem to $H_c^*$ and $g_c^*$ one obtains
\begin{eqnarray*}
&&\rho((y,g,H,w),(y_c^*,g_c^*,H_c^*,w_c^*))\\
&\leq& \epsilon +
\max_{l=1,\dots,N_{B}}\sup_{\|\bz-\bz_l\|\leq \epsilon}
|H(y+g(\bz))-H_c^*(y_c^*+g_c^*(\bz))|\\
&\leq& \epsilon +\|H-H_c^*\|_\infty+\sup_{x}h_c^*(x)\|g-g_c^*\|_\infty\\
&&{}+\max_{l=1,\dots,N_{B}}\sup_{\|\bz-\bz_l\|\leq \epsilon}
|H_c^*(y+g_c^*(\bz))-H_c^*(y_c^*+g_c^*(\bz))|\\
&\leq& (2+(1+\eta)C)\epsilon+\max_{l=1,\dots,N_{B}}
|H_c^*(y_{l,c}^{*2}+g_c^*(\bz_l))-H_c^*(y_{l,c}^{*1}+g_c^*(\bz_l))|\\
&\leq& (4+(1+\eta)C)\epsilon. 
\end{eqnarray*}
For the sake of brevity we omit the proof of the remaining (simpler) conditions for application of Theorem 2.11.9 in Van der Vaart and Wellner (1996). 
An application of the theorem gives asymptotic $\rho$-equicontinuity of the process $E_n$ which implies that $\sup_{y,g,H}|E_n(y,g,H,w) - E_n(y,g,H,0)| \rightarrow 0$ in probability as $w \rightarrow 0$.
\hfill $\Box$

\begin{Lemma} \label{lem2-linear}
Assume the linear model (\ref{mod}) and \ref{AL.1}, \ref{AL.2}, \ref{CL.1} and \ref{CL.2}.   Then, conditionally on the data, the process $T_{n2}(y)$, $y \in \R$, defined in (\ref{T_n2-lin}), converges weakly to the Gaussian process $W(y)$ in probability, where $W(y)$ is  defined in (\ref{covW-lin}).
\end{Lemma}

\noindent
{\bf Proof.} Note that we have
$$ T_{n2}(y)=n^{-1/2} \sumi \big(I\{\hat\eps_{s,0,i}^* \le y\} - \hat F_{s,n}(y) \big)$$ 
with $\hat\eps_{s,0,i}^*=\eps_{s,i}^*+\bx_{ni}^T(\hat\bbeta_n-\hat\bbeta_{0,n}^*)$ with $\hat\bbeta^*_{0,n}$ from (\ref{hat-beta*}). The process $W_n^*$ considered by Koul and Lahiri (1994) (in the case of least squares estimators, i.e.\ $\psi(x)=x$, and with weights $b_{ni}=n^{-1/2}$) corresponds to $T_{n2}\circ \hat F_{s,n}^{-1}$, but replacing $\hat\eps_{s,0,i}^*$ by $\hat \eps_{s,i}^*= \eps_{s,i}^*+\bx_{ni}^T(\hat\bbeta_n-\hat\bbeta_{s,n}^*)$ with 
$\hat\bbeta_{s,n}^*=\hat\bbeta_n+(\bX_n^T\bX_n)^{-1}\bX_n^T\beps_{s,n}^*$ (based on the smooth bootstrap residuals). We have to show that the use of the different parameter estimator $\hat\bbeta_{0,n}^*$ (instead of $\hat\bbeta_{s,n}^*$) does not change the result of Theorem 2.1 in Koul and Lahiri (1994). Assumptions (A.1), (A.2) and (A.3)(iii) in the aforementioned paper are fulfilled under our assumptions. Assumption (A.3)(ii) in our case reads as 
$$E^*[(\eps_{in}^*)^k]\longrightarrow E[(\eps_{11})^k] \mbox{ for }n\to\infty \mbox{ a.s.},\quad k=1,2. $$
The condition is fulfilled for the smooth bootstrap errors $\eps_{in}^*=\eps_{s,i}^*$ by Proposition 2.1 in Koul and Lahiri (1994). Further the same condition for the non-smooth bootstrap errors $\eps_{in}^*=\eps_{0,i}^*$ is equivalent to 
\begin{equation}\label{moments}
\frac1n\sumi (\hat\eps_{in})^k\longrightarrow E[(\eps_{11})^k] \mbox{ for }n\to\infty \mbox{ a.s.},\quad k=1,2,
\end{equation}
which is a standard result in least squares estimation. 
Assumption (A.3)(iv) in Koul and Lahiri (1994) for $\hat\bbeta^*_{0,n}$ is valid by (\ref{hat-beta*}), but gives an expansion in terms of $\eps_{0,i}^*$ (as opposed to an expansion in the smooth $\eps_{s,i}^*$). Lemma 3.1 in Koul and Lahiri (1994) (with $Z_{ni}=\eps_{s,i}^*$, $H_n=\hat F_{s,n}$, $b_{ni}=n^{-1/2}$, $d_{ni}^T=\bx_{ni}^T(\bX_n^T\bX_n)^{-1/2}$, $u=(\bX_n^T\bX_n)^{1/2}(\hat\bbeta^*_{0,n}-\hat\bbeta_n)$) gives the representation 
\begin{eqnarray*}
T_{n2}(\hat F_{s,n}^{-1}(t)) &=& n^{-1/2}\sumi (I\{\eps_{s,i}^*\leq \hat F_{s,n}^{-1}(t)\}-t)\\
&&{} + \hat f_{s,n}(\hat F_{s,n}^{-1}(t))n^{-1/2}\sumi x_{ni}^T (\hat\bbeta_{0,n}^*-\hat\bbeta_n) +o_P(1)
\end{eqnarray*}
uniformly in $t\in [0,1]$, where $\hat f_{s,n}$ is the density corresponding to $\hat F_{s,n}$.  From (\ref{hat-beta*}) it now follows that
\begin{eqnarray*}
T_{n2}(\hat F_{s,n}^{-1}(t)) &=& n^{-1/2}\sumi (I\{\eps_{s,i}^*\leq \hat F_{s,n}^{-1}(t)\}-t)+ \hat f_{s,n}(\hat F_{s,n}^{-1}(t))n^{-1/2}\sum_{j=1}^n \tilde x_{nj} \eps_{0,j}^* +o_P(1)
\end{eqnarray*}
uniformly in $t\in [0,1]$ with 
\begin{eqnarray}\label{tilde_x}
\tilde x_{nj}&=& \sum_{i=1}^n\bx_{ni}^T(\bX_n^T\bX_n)^{-1}\bx_{nj}.
\end{eqnarray} 
Note  that assumptions \ref{2} and \ref{3} imply 
\begin{eqnarray} \label{rem1} 
\frac{1}{n}\sum_{j=1}^n\tilde x_{nj}^2=\frac{1}{n}\sum_{j=1}^n\tilde x_{nj}&\longrightarrow& \bold{m}^T\boldsymbol{\Sigma}^{-1}\bold{m}\quad\mbox{for $n\to\infty$}.
\end{eqnarray} 
Because of the continuity of the distribution $\hat F_{s,n}$ with support $\R$ (see condition \ref{CL.1}) one obtains directly $T_{n2}(y) = \tilde T_{n2}(y) + o_P(1)$ uniformly in $y\in\R$, where 
$$ \tilde T_{n2}(y) = n^{-1/2} \sumi \big(I\{\eps_{s,i}^* \le y\} - \hat F_{s,n}(y) + \hat f_{s,n}(y) \tilde x_{ni}\eps_{0,i}^*\big). $$
To show conditional weak convergence of $\tilde T_{n2}$ to $W$ in probability one can follow the steps of the proof of Theorem 2 in Neumeyer (2009). To this end one needs  uniform convergence of $\hat f_{s,n}$ to $f$, which is assumption (A.3)(iii) in Koul and Lahiri (1994) and follows from their Proposition 2.1 under our assumption (see also the proof of Lemma \ref{lem1-linear}). 
Further one needs to apply (\ref{moments}) and (\ref{rem1}) and to show convergence of $E^*[I\{\eps_{s,i}^* \le y\}\eps_{0,i}^*]$ to $E[I\{\eps\le y\}\eps]$ in probability, which is analogous to the last part of the proof of Lemma \ref{lem2}.
\hfill $\Box$

\begin{Lemma} \label{lem3-linear}
Assume the linear model (\ref{mod}) and \ref{AL.1}, \ref{AL.2}, \ref{CL.1} and \ref{CL.2}.    Then, conditionally on the data, $\sup_y |T_{n3}(y)| = o_P(1)$, in probability, where
 $$T_{n3}(y) =n^{1/2} \big(\hat F_{s,n}(y) - \hat F_{0,n}(y)\big) .$$
\end{Lemma}

\noindent
{\bf Proof.} We will use the expansion given in 
 Theorem 6.2.1 in Koul (2002, p.\ 232), i.e.
\begin{eqnarray*}
\hat F_{0,n}(y) &=& \frac{1}{n}\sumi I\{\eps_{ni} \le y\}+f(y)\frac{1}{n}\sum_{i=1}^n\bx_{ni}^T(\hat{\bbeta}_n-\bbeta)+o_P(n^{-1/2})\\
&=& \frac{1}{n}\sumi I\{\eps_{ni} \le y\}+f(y)\frac{1}{n}\sum_{j=1}^n
\tilde x_{nj}\eps_{nj}+o_P(n^{-1/2}),
\end{eqnarray*}
uniformly with respect to $y\in\R$. Here the second equality follows from (\ref{M}) and $\tilde x_{nj}$ was defined in (\ref{tilde_x}). 
Recall that $\hat F_{s,n}(y) =\int \hat F_{0,n}(y-vs_n)\,dL(v) $. 
With (\ref{rem1}) we have
$$ \frac{1}{n}\sum_{j=1}^n\tilde x_{nj}\eps_{nj} \int \big[f(y-vs_n) - f(y)\big] \, dL(v)= O_P(n^{-1/2})\int \big[f(y-vs_n) - f(y)\big] \, dL(v) =o_P(n^{-1/2})$$
by the properties of the density $f$ and the kernel $\ell$. 
It remains to show that
\begin{eqnarray*}
\int \frac{1}{n}\sumi I\{\eps_{ni} \le y+s_nv\} \, dL(v) &=&\frac{1}{n}\sumi I\{\eps_{ni} \le y\} + o_P(n^{-1/2}) 
\end{eqnarray*}
uniformly in $y$. 
Note that  the term on the left hand side of the equation is the empirical distribution function of iid data $\eps_{ni}$, $i=1,\dots,n$, smoothed with kernel $\ell$ and bandwidth $s_n$, whereas the term on the right hand side is the classical edf of those data.
The theorem in Van der Vaart (1994) gives the desired  $o_P(n^{-1/2})$-rate of their difference. To see this consider the discussion on p.\ 502 of the aforementioned paper for $\sigma_n=s_n=o(n^{-1/(2k)})$ for $k=2$, $\mu_n$ equal to the distribution of $s_nU$, where $U$ has density $\ell$, and for the function class $\mathcal{F}=\{\eps\mapsto f_z(\eps)= I\{\eps\leq z \}\mid z\in\R\}$, which results in twice differentiable functions $y\mapsto Pf_z(\varepsilon+y)=F(z-y)$ with bounded derivatives. 
\hfill $\Box$

\bigskip

\newpage

\begin{center}
{\Large \bf References}
\end{center}

\bib Akritas, M. G. and Van Keilegom, I. (2001).
Nonparametric estimation of the residual distribution.
{\it Scand. J. Statist.}  {\bf 28}, 549--567.

\bib Azzalini, A. (1985). 
A class of distributions which includes the normal ones. 
{\it Scand. J. Statist.} {\bf 12}, 171--178.

\bib Cheng, G. and Huang, J.Z. (2010). 
Bootstrap consistency for general semiparametric M-estimation. 
{\it Ann. Statist.} {\bf 38}, 2884--2915.


\bib De Angelis, D., Hall, P. and Young, G. A. (1993). Analytical and bootstrap approximations to estimator
distributions in $L^1$ regression. {\it J. Amer. Statist. Assoc.} {\bf 88}, 1310--1316.

\bib Durbin, J. (1973). 
Weak convergence of the sample distribution function when parameters are estimated. 
{\it Ann. Statist.} {\bf 1}, 279--290.

\bib Einmahl, J. and Van Keilegom, I. (2008). 
Specification tests in nonparametric regression. 
{\it J. Econometrics} {\bf 143}, 88--102.

\bib  Hall, P., DiCiccio, T.J. and Romano, J.P. (1989).
On smoothing and the bootstrap
{\it Ann. Statist.} {\bf 17}, 692--704. 

\bib Huskova, M. and Meintanis, S.G. (2009).
Goodness-of-fit tests for parametric regression models based on empirical
characteristic functions. 
{\it Kybernetika} {\bf 45}, 960--971.

\bib Koul, H.L. (1987). Tests of goodness-of-fit in linear regression. 
{\it Goodness-of-fit (Debrecen, 1984) Colloq. Math. Soc. Janos Bolyai} {\bf 45},  
279--315, North-Holland, Amsterdam.

\bib Koul, H.L. (2002). 
{\it Weighted Empirical Processes in Dynamic Nonlinear Models (Second Edition).} Springer, New York.

\bib Koul, H.L. and Lahiri, S.N. (1994).
On bootstrapping M-estimated residual processes in multiple linear regression models.
{\it J. Multiv. Anal.} {\bf 49}, 255--265.

\bib Loynes, R.M. (1980). 
The empirical sample distribution function of residuals from generalized regression. 
{\it Ann. Statist.} {\bf 8}, 285--298.

\bib Mora, J. (2005). Comparing distribution functions of errors in linear models: a nonparametric approach. {\it Statist. Probab. Lett.} {\bf 73}, 425--432. 

\bib M\"uller, U.U., Schick, A. and Wefelmeyer, W. (2004). 
Estimating linear functionals of the error distribution in nonparametric regression.
{\it J. Statist. Planning Inf.} {\bf  119}, 75--93. 

\bib Neumeyer, N.\ \rm (2006). Bootstrap procedures for empirical processes of nonparametric residuals. Habilitationsschrift, Ruhr-Universit\"at Bochum.\\ available at http://www.math.uni-hamburg.de/home/neumeyer/habil.ps 

\bib Neumeyer, N. (2008).
A bootstrap version of the residual-based smooth empirical distribution function.
{\it J. Nonparametr. Statist.} {\bf 20}, 153--174.

\bib Neumeyer, N. (2009). 
Smooth residual bootstrap for empirical processes of nonparametric regression residuals. 
{\it Scand. J. Statist.}  {\bf 36}, 204--228.


\bib Neumeyer, N., Dette, H. and Nagel, E.-R. (2005). 
A note on testing symmetry of the error distribution in linear regression models. {\it J. Nonparametr. Stat.} {\bf  17}, 697--715. 

\bib Neumeyer, N. and Dette, H. (2007). 
Testing for symmetric error distribution in nonparametric regression models. 
{\it Statist. Sinica} {\bf 17}, 775--795.


\bib Pardo-Fern\'andez, J.C., Van Keilegom, I. and Gonz\'alez-Manteiga, W. (2007). 
Testing for the equality of $k$ regression curves. 
{\it Statist. Sinica} {\bf 17}, 1115--1137.

\bib Racine, J. and Van Keilegom, I. (2017). 
A smooth nonparametric, multivariate, mixed-data location-scale test 
(submitted). 

\bib  Silverman, B.W.  and Young, G.A. (1987). The bootstrap: to smooth or not to smooth?
{\it Biometrika} {\bf 74}, 469--479. 


\bib Van der Vaart, A. W. (1994). Weak convergence of smoothed empirical processes. {\it  Scand. J. Statist.} {\bf  21}, 501--504. 

\bib Van der Vaart, A. W. and Wellner, J. A. (1996).
\textit{Weak Convergence and Empirical Processes}.
Springer-Verlag, New York.

\bib Van Keilegom, I. and Akritas, M.G. (1999). 
Transfer of tail information in censored regression models. 
{\it Ann. Statist.} {\bf 27}, 1745--1784.

\bib Van Keilegom, I., Gonz\'alez-Manteiga, W. and S\'anchez-Sellero, C. (2008). 
Goodness-of-fit tests in parametric regression based on the estimation of the error distribution. 
{\it TEST} {\bf 17}, 401--415.

\end{document}